\documentclass[11pt]{amsart}
\usepackage{amssymb}
\usepackage[sectiontheorem]{basicmacros}

\newcommand{\nk}[1]{\ensuremath{#1_{n}^{k}}}
\newcommand{\ank}{\nk{a}}
\newcommand{\bnk}{\nk{b}}
\newcommand{\intnk}{\opint{\ank}{\bnk}}

\newcommand{\fII}{\selfmap{f}{\Interval}}

\newcommand{\fold}[1]{\ensuremath{#1}-fold}
\newcommand{\mfold}{\fold{m}}
\newcommand{\comfold}{cocountably \mfold}
\newcommand{\msec}{$m$-section}

\newcommand{\psijof}[1]{\psisof{j}{#1}}
\newcommand{\dist}{\ensuremath{\operatorname{dist}}}
\newcommand{\distof}[2]{\ensuremath{\dist(#1,#2)}}
\newcommand{\mdistof}[2]{\ensuremath{\operatorname{mdist}(#1,#2)}}
\newcommand{\varof}[2]{\ensuremath{\norm{#1}_{#2}}}

\renewcommand{\ms}[1]{\ensuremath{m_{#1}}}
\newcommand{\my}{\ms{y}}
\newcommand{\mys}[1]{\ms{\ys{#1}}}
\newcommand{\myi}{\mys{i}}
\renewcommand{\Ms}[1]{\ensuremath{M_{#1}}}
\newcommand{\My}{\Ms{y}}
\newcommand{\Mys}[1]{\Ms{\ys{#1}}}
\newcommand{\Myi}{\Mys{i}}

\newcommand{\ent}[1]{\ensuremath{h_{top}(#1)}}
\newcommand{\entof}[1]{\ensuremath{\operatorname{ent}(#1)}}
\newcommand{\vent}[1]{\entof{\itinset{#1}}}

\newcommand{\tJs}[1]{\ensuremath{\tilde{J}_{#1}}}

\newcommand{\cocount}{cocountable}
\newcommand{\leftreg}[1]{left $#1$-regular }
\newcommand{\rightreg}[1]{right $#1$-regular }
\newcommand{\reg}[1]{$#1$-regular}
\newcommand{\Regset}{\ensuremath{\mathfrak{C}}}
\newcommand{\Lreg}[1]{\ensuremath{\Regset_{#1}(f,\ell)}}
\newcommand{\Lregm}{\Lreg{m}}
\newcommand{\Rreg}[1]{\ensuremath{\Regset_{#1}(f,r)}}
\newcommand{\Rregm}{\Rreg{m}}
\newcommand{\Reg}[1]{\ensuremath{\Regset_{#1}(f)}}
\newcommand{\Regm}{\Reg{m}}

\newcommand{\Shiftsystem}{\ensuremath{\mathcal{H}}}
\newcommand{\Shiftpiece}{\ensuremath{\mathbb{H}}}
\newcommand{\Shiftsubset}{\ensuremath{H}}
\newcommand{\Shiftset}[1]{\ensuremath{\Shiftsubset_{#1}}}

\newcommand{\Domain}{\ensuremath{\mathbb{D}}}

\newcommand{\DomH}{\Domain}
\newcommand{\DomsH}[1]{\ensuremath{\Domain_{#1}}}

\newcommand{\inDomsH}[2]{\inside{#2}{\DomsH{#1}}}

\newcommand{\kerof}[1]{\ensuremath{\mathfrak{K}(#1)}}
\newcommand{\cZ}{\ensuremath{\mathfrak{Z}}}
\newcommand{\cent}[1]{\ensuremath{\cZ(#1)}}
\newcommand{\kerH}{\kerof{\Shiftsystem}}
\newcommand{\centH}{\cent{\Shiftsystem}}
\newcommand{\centOH}{\ensuremath{\cZ_{0}(\Shiftsystem)}}

\newcommand{\shiftspace}{\ensuremath{\Omega}}
\newcommand{\mshiftspace}{\ensuremath{\shiftspace_{m}}}
\newcommand{\cylinders}[1]{\ensuremath{\mshiftspace(#1)}}

\newcommand{\code}{\ensuremath{a}}

\newcommand{\ads}[1]{\ensuremath{\code_{#1}}}
\newcommand{\adset}[1]{\ensuremath{\alpha(#1)}}
\newcommand{\mult}[1]{\ensuremath{\eta(#1)}}
\newcommand{\itinset}[1]{\ensuremath{\shiftspace(#1)}}
\newcommand{\nitinset}[2]{\ensuremath{\shiftspace(#1)(#2)}}
\newcommand{\realizeset}[1]{\ensuremath{\Pi(#1)}}

\newcommand{\iss}[2]{\ensuremath{i_{#1}^{(#2)}}}
\newcommand{\jss}[2]{\ensuremath{j_{#1}^{(#2)}}}

\newcommand{\ijs}[1]{\ensuremath{(\is{#1},\js{#1})}}

\newcommand{\T}{\ensuremath{T}}
\newcommand{\Tpr}{\ensuremath{\T\pr}}
\newcommand{\fTT}{\selfmap{f}{\T}}
\newcommand{\Peri}{\ensuremath{Per(f)}}
\newcommand{\Fix}{\ensuremath{Fix(f)}}
\newcommand{\valence}{\ensuremath{\nu}}
\newcommand{\inT}[1]{\ensuremath{#1\in\T}}

\newcommand{\bpt}{branchpoint}
\newcommand{\Branchpts}{\ensuremath{\mathsf{B}(\T)}}
\newcommand{\Vertices}{\ensuremath{\mathsf{V}(\T)}}

\newcommand{\y}{\ensuremath{y}}
\newcommand{\vBranchsymbol}{\ensuremath{B}}
\newcommand{\Branch}[1]{\ensuremath{\vBranchsymbol_{#1}} }
\newcommand{\Branchi}{\Branch{i}}
\newcommand{\Branchat}[2]{\ensuremath{\vBranchsymbol_{#1}(#2)}}
\newcommand{\yBranch}[1]{\ensuremath{\Branchat{#1}{y}}}
\newcommand{\U}{\ensuremath{U}}
\newcommand{\Ubranch}[1]{\ensuremath{\U_{#1}}}
\newcommand{\Usi}{\Ubranch{i}}

\newcommand{\tU}{\ensuremath{\tilde{\U}}}
\newcommand{\tUbranch}[1]{\ensuremath{\tU_{#1}}}
\newcommand{\tUsi}{\tUbranch{i}}

\newcommand{\Vof}[1]{\ensuremath{V(#1)}}

\newcommand{\upbranchsymbol}{\ensuremath{\mathcal{B}}}
\newcommand{\Upbranches}[2]{\ensuremath{\upbranchsymbol^{+}_{#2}(#1)}}
\newcommand{\Downbranches}[2]{\ensuremath{\upbranchsymbol^{-}_{#2}(#1)}}

\newcommand{\colorit}{color itinerary}
\newcommand{\colorits}{color itineraries}
\newcommand{\colorword}{\ensuremath{\mathfrak{c}}}
\newcommand{\cword}[1]{\ensuremath{\colorword(#1)}}
\newcommand{\Numcolors}[1]{\ensuremath{\mathcal{N}_{#1}}}
\newcommand{\Numn}{\Numcolors{n}}
\newcommand{\Graph}[2]{\ensuremath{\Gamma(#1,#2)}}
\newcommand{\bcgraph}{\Graph{b}{c}}
\newcommand{\bcs}[1]{\ensuremath{(b,c)_{#1}}}
\newcommand{\bcn}{\bcs{n}}
\newcommand{\gamp}{\ensuremath{\gamma\pr}}
\newcommand{\gammin}{\gams{-}}
\newcommand{\Pname}{\ensuremath{\mathcal{C}}}
\newcommand{\Pnum}[3]{\ensuremath{\Pname_{#1}(#2,#3)}}
\newcommand{\Pnumbc}[1]{\Pnum{#1}{b}{c}}
\newcommand{\Pnumbcn}{\Pnum{n}{b}{c}}

\newcommand{\Mnum}[3]{\ensuremath{\mathcal{M}_{#1}(#2,#3)}}
\newcommand{\Mbcnum}[1]{\Mnum{#1}{b}{c}}
\newcommand{\Mbcnump}{\Mbcnum{p}}

\newcommand{\Snum}[3]{\ensuremath{\mathcal{S}_{#1}(#2,#3)}}
\newcommand{\Snumbc}[1]{\Snum{#1}{b}{c}}
\newcommand{\Snumbcp}{\Snumbc{p}}

\newcommand{\Percor}{\ensuremath{\mathcal{P}}}

\newcommand{\len}[1]{\ensuremath{\mathrm{len}(#1)}}
\newcommand{\lenlam}{\ensuremath{\len{\lambda}}}
\newcommand{\qfactor}{\ensuremath{(1-\frac{1}{q})}}

\newcommand{\bgraph}{\ensuremath{\mathcal{B}}}
\newcommand{\bcovers}[2]{\covers{\bgraph}{#1}{#2}}
\newcommand{\cover}[4]{\ensuremath{(#1,#2)\rightarrow( #3,#4)}}

\newcommand{\admitless}{\ensuremath{\prec}}
\newcommand{\admitmore}{\ensuremath{\succ}}
\newcommand{\concat}[2]{\ensuremath{{#1}{#2}}}

\newcommand{\Ph}{\ensuremath{\Phi}}

\newcommand{\covers}[3]{\ensuremath{#2\overset{#1}{\rightarrow}{#3}}}

\newcommand{\Hs}[1]{\ensuremath{H_{#1}}}

\newcommand{\Hstar}{\Hs{*}}

\newcommand{\Hone}{\Hs{1}}
\newcommand{\Hm}{\Hs{m}}
\newcommand{\Hboth}{\Hs{\#}}
\newcommand{\Hplus}{\Hs{+}}
\newcommand{\Hminus}{\Hs{-}}

\newcommand{\Gsys}{\ensuremath{\mathfrak{G}}}
\newcommand{\Gset}{\ensuremath{G}}
\newcommand{\Gss}[2]{\ensuremath{\Gset_{#1}^{#2}}}
\newcommand{\Gij}{\Gss{i}{j}}

\renewcommand{\ip}{\ensuremath{i\pr}}
\newcommand{\jp}{\ensuremath{j\pr}}

\newcommand{\pp}{\ensuremath{p\pr}}

\newcommand{\ps}[1]{\ensuremath{p_{#1}}}
\newcommand{\qps}[1]{\ensuremath{q\pr_{#1}}}
\newcommand{\qpp}{\ensuremath{q^{\prime\prime}}}
\newcommand{\qpps}[1]{\ensuremath{\qpp_{#1}}}

\newcommand{\edgeitin}[1]{\ensuremath{\mathsf{ei}(#1)}}
\newcommand{\Es}[1]{\ensuremath{E_{#1}}}
\newcommand{\Eps}[1]{\ensuremath{E\pr_{#1}}}
\newcommand{\edgequiv}[2]{\ensuremath{#1\sim #2}}
\newcommand{\nedgequiv}[2]{\ensuremath{#1\not\sim #2}}
\newcommand{\edgeclass}[1]{\ensuremath{[#1]}}
\newcommand{\edgehull}[1]{\hull{#1}}
\newcommand{\edgehullorbsymbol}{\ensuremath{Z}}
\newcommand{\edgehullorb}[1]{\ensuremath{\edgehullorbsymbol(#1)}}
\newcommand{\sepsetsymbol}{\ensuremath{W}}
\newcommand{\sepset}[1]{\ensuremath{\sepsetsymbol(#1)}}
\newcommand{\Sepset}{\ensuremath{W}}

\newcommand{\Psp}[1]{\ensuremath{P_{#1}(p)}}
\newcommand{\Pspmin}{\ensuremath{P^{-}_{0}(p)}}
\newcommand{\Pspmax}{\ensuremath{P^{+}_{0}(p)}}
\newcommand{\Pspmaxpo}{\ensuremath{P^{+}_{0}(\ps{0})}}
\newcommand{\Pspeither}[1]{\ensuremath{P^{\pm}_{#1}(p)}}
\newcommand{\Zp}{\ensuremath{\edgehullorb{p}}}
\newcommand{\Zpo}{\ensuremath{\edgehullorb{\psub{0}}}}
\newcommand{\Zq}{\ensuremath{\edgehullorb{q}}}
\newcommand{\Zs}[1]{\ensuremath{\edgehullorbsymbol_{#1}(p)}}
\newcommand{\Zcompl}{\ensuremath{\T\setminus \Zp}}
\newcommand{\Zcomplpo}{\ensuremath{\T\setminus \Zpo}}
\newcommand{\Zcomplq}{\ensuremath{\T\setminus \Zq}}
\newcommand{\Zis}[1]{\Zs{\is{#1}}}

\newcommand{\Centcomp}{\ensuremath{C(p)}}
\newcommand{\Percomp}{\ensuremath{P(p)}}

\newcommand{\zpls}[1]{\ensuremath{z^{+}_{#1}}}
\newcommand{\zmins}[1]{\ensuremath{z^{-}_{#1}}}
\newcommand{\zpms}[1]{\ensuremath{z^{\pm}_{#1}}}

\newcommand{\Usymbol}{\ensuremath{U}}
\newcommand{\Uss}[2]{\ensuremath{\Usymbol^{#1}_{#2}(p)}}
\newcommand{\Ups}[1]{\Uss{+}{#1}}
\newcommand{\Ums}[1]{\Uss{-}{#1}}
\newcommand{\Upms}[1]{\Uss{\pm}{#1}}

\newcommand{\hull}[1]{\ensuremath{\langle #1 \rangle}}

\newcommand{\zpp}{\ensuremath{z^{\prime\prime}}}

\newcommand{\Ftil}[1]{\ensuremath{\tilde{F}_{#1}}}

\title{Topological Entropy of \mfold{} Maps on Trees}

\author{ Jozef Bobok}
	\thanks{First author supported  by 
	MYES of the Czech Republic via Contract MSM 684 077 0010.}
	\thanks{Both authors thank the Departament de
	Matem$\grave{\text{a}}$tiques,  Universitat
	Aut$\grave{\text{o}}$noma de Barcelona for its kind hospitality and
	partial support during separate visits during 2003-4.}
	\address{KM FSv. \v{C}VUT, Th\'akurova 7, 166 29 Praha 6,
        Czech Republic}
        \email{bobok@mat.fsv.cvut.cz}
\author{Zbigniew Nitecki}
	\address{Department of Mathematics, 
        Tufts University, Medford, MA 02155, USA}
        \email{zbigniew.nitecki@tufts.edu}
\date{\today}

\begin{document}

\begin{abstract}
	We establish the analogue for maps on trees of the result established in 
	\cite{Bobok:2fold, Bobok:mfold} for interval maps, that a continuous self-map for which all but
	countably many points have at least $m$ preimages (and none have less than two)
	has topological entropy bounded below by $\log m$.  
\end{abstract} 

\maketitle

\section{Introduction}
 This paper continues the investigation in \cite{Bobok:2fold, Bobok:mfold, BobokZN:mfold} of the relation between the topological entropy \ent{f} of a continuous map \fXX{} of a compact space to itself and the number of preimages of points under $f$.
 
 We thank Lluis Alsed\`a for extensive and helpful discussions during the writing of this paper.
 
 A continuous map \fXX{} is \deffont{$\boldsymbol{m}$-fold} on a subset $\Y\subset\X$ if every point in \Y{} has at least $m$ preimages in \X{}.  We drop the reference to the subset  \Y{} when it is the whole space. 
 
An argument due to Misiurewicz and Przytycki \cite{MisPrzy:degree} shows that any $C^{1}$ map on any compact manifold which is \mfold{} on the set of its regular (non-critical) values satisfies the estimate
		\begin{equation}\label{eqn:main}
			\ent{f}\geq\log m.
		\end{equation}
This argument is detailed in \cite{BobokZN:mfold}.  Ethan Coven \cite{Coven:conj} conjectured that the differentiability condition could be replaced by continuity for a \fold{2} map of the interval;  this was established in \cite{Bobok:2fold}.  

A simple example \cite{BobokZN:mfold} shows that for a continuous map of the interval, failure of the \fold{2} condition at a single point can allow entropy zero; the same example can be adapted to create an \mfold{} map on the circle with zero entropy.  However, in \cite{Bobok:mfold} it was shown that for a map of the interval, once the \fold{2} condition is assumed to hold everywhere, the estimate \eqref{eqn:main} is guaranteed with any higher value of $m$ as soon as the \mfold{} condition holds on the complement of a countable set.  We call a continuous map \fXX{} \deffont{cocountably $\boldsymbol{m}$-fold} if it is \fold{2} on \X{} and \mfold{} on a cocountable set $\Y\subset\X$ (\ie{} the complement $\X\setminus\Y$ has at most countably many points). In \cite{BobokZN:mfold} it was shown that \refer{eqn}{main} holds for a \comfold{} self-map of the circle provided there is a positive lower bound on the diameter of all preimage sets of points.

In this paper, we establish the following extension to trees of the result established by the second author \cite{Bobok:mfold} for maps of the interval:
 \begin{theorem}\label{thm:main} Suppose \fTT{} is a \comfold{} map of a (finite) tree to itself.  	Then \refer{eqn}{main} holds:
 	$$\ent{f}\geq\log m.$$
\end{theorem}
 
Our strategy follows that of \cite{Bobok:2fold, Bobok:mfold}, with some modifications to adapt it to more general trees.  In \refer{sec}{EntViaShift} we describe a modification of \cite[Theorem 4.8]{BobokZN:mfold}, which systematizes the strategy of \cite{Bobok:2fold, Bobok:mfold} as an abstract set of conditions on a kind of weak ``horseshoe'' (which we refer to as a \emph{shift system}) that guarantee the estimate \eqref{eqn:main}. This is general symbolic dynamics, and is not specific to tree maps. Then in \refer{sec}{background} we formulate a scheme for linearly ordering the points of a tree, which takes the place of the linear ordering on an interval in our arguments. The hypotheses of the symbolic result are established for any \mfold{} map on a tree in \S\S\ref{sec:kernel}-\ref{sec:center}, and everything is combined to prove our main result in \refer{sec}{main}.

We  close this section with an example which shows that the set of points where the \mfold{} condition fails in the hypotheses of \refer{thm}{main} cannot be allowed to be uncountable, even if it is nowhere dense (at least not without further conditions).

\begin{prop}\label{prop:example}
	For each integer $m>0$ there exists a map \fII{} of the 
	interval $\Interval=\clint{0}{1}$ to itself such that
	\begin{enumerate}
		\item $f$ is globally \fold{2};
		\item $f$ is \mfold{} on a set $\Y=I\setminus K$, where $K$ is a nowhere dense, closed
		(uncountable) set;
		\item $\ent{f}=\log 2$.
	\end{enumerate}
\end{prop}

\begin{proof}
	We begin with a ``flattened tent map'' \selfmap{g}{\Interval}, taking both endpoints to $0$, 
	taking a central interval \clint{a}{b} to the right endpoint, and affine (or even just strictly monotone)
	 on each of the 
	complementary intervals.  For example, taking $a=\frac{1}{3}$, 
	$b=\frac{2}{3}$, we
	can define $g$ by
	\begin{equation*}
		g(x)\eqdef
		\begin{cases}
			3x  & \text{for }0\leq x\leq \frac{1}{3}\\
			1      & \text{for  }\frac{1}{3}\leq x\leq\frac{2}{3} \\
			3(1-x)      & \text{for }\frac{2}{3}\leq x\leq 1.
		\end{cases}
	\end{equation*}
	Clearly, $\ent{g}=\log 2$.
	
	Now, the iterated preimages of the central interval \opint{a}{b} consist of disjoint open intervals
		$$\gpre{k}{\opint{a}{b}}=\bigcup_{n=1}^{2^{k}} \intnk{},
		\quad \inflist{k}{0,1,2}$$
	whose union is a dense open set
		$$\Y=\bigcup_{k=0}^{\infty}\pre{g}{k}{\opint{a}{b}};$$
	the complement of \Y{} is a Cantor-like set $K\subset\Interval$ 
	(with $a=\frac{1}{3}$ and $b=\frac{2}{3}$, $K$ is the classical middle-third Cantor set).
	We modify $g$ on each of the intervals \intnk{}, $k>0$, to a map \fII{} mapping \intnk{} 
	and its mirror image, \opint{1-\bnk}{1-\ank}, onto
	\gof{\intnk} in an $m$-fold manner.
	
	Then clearly $f$ is \mfold{} on the set \Y{}, but since $\ftoof{k+1}{\intnk}=\single{1}$, 
	each of the intervals \intnk{} is wandering;  thus 
		$$\ent{f}= \ent{f|K}=\ent{g|K}=\log 2.$$
\end{proof}


\section{Entropy via Shift Systems}\label{sec:EntViaShift}
The material of this section is not specific to trees, and closely follows \cite[\S\S3-4]{BobokZN:mfold}. 
However, a modification of the definition of a ``locally dividing'' set was needed for our purposes here.  We shall sketch many of the arguments, referring the reader to the exposition in \cite{BobokZN:mfold} for details, but provide more detailed proofs where the modification mentioned above requires them.

An \deffont{$\boldsymbol{m}$-shift system} for a map \fXX{} is a collection
 \begin{equation*}
	\Shiftsystem{}=\single{\Shiftset{1},...,\Shiftset{m}}
\end{equation*} 
of $m$ 
nonempty (but not necessarily closed or disjoint) sets $\Shiftset{i}\subset\X$ satisfying 
\begin{equation}\label{eqn:shiftsystem}
	\fof{\Shiftset{i}}\supseteq\Shiftpiece\eqdef\Shiftset{1}\cup\dots\cup\Shiftset{m}
		\quad\text{for }\ilist{i}{1}{m}.
\end{equation}

The \deffont{address set} of \inX{x} is
\begin{equation*}
	\adset{x}\eqdef\setbld{\code\in\setlist{1}{m}}{x\in\Shiftset{\code}}
\end{equation*}
and its cardinality $\mult{x}$ is the \deffont{multiplicity} of \Shiftsystem{} at $x$.  The set of points
with positive multiplicity is precisely \Shiftpiece; 
we define the \deffont{kernel} \resp{\deffont{center}} of
\Shiftsystem{} to be the set of points with multiplicity greater than one \resp{equal to $m$}
\begin{eqnarray*}
	\kerof{\Shiftsystem}\eqdef\setbld{x}{\mult{x}>1}&=&\bigcup_{i\neq j}\Shiftset{i}\cap\Shiftset{j}\\
	\cent{\Shiftsystem}\eqdef\setbld{x}{\mult{x}=m}&=&\bigcap_{i=1}^{m}\Shiftset{i}.
\end{eqnarray*}
We also define the \deffont{core} of \Shiftsystem{} as the set of points whose orbit remains in \centH{} for all time:
\begin{equation*}
	\centOH\eqdef\setbld{\inX{x}}{\ftoof{k}{x}\in\centH\text{ for all }k}
		=\bigcap_{k=0}^{\infty}\fpre{k}{\centH}.
\end{equation*}
Obviously, $\centOH\subseteq\centH\subseteq\kerH\subseteq\Shiftpiece$.  We will call the 
shift system \Shiftsystem{} \deffont{nontrivial} if $\Shiftset{i}\setminus\centOH\neq\emptyset$ for all $i$, and \deffont{closed} if each \Shiftset{i} is a closed subset of \X{}.

A closed $m$-shift system with \emph{empty} kernel ($\kerH=\emptyset$) is usually called a \deffont{horseshoe} in the context of maps of the interval \cite{AlsedaLlibreMis};  the \emph{itinerary} of a point with respect to a horseshoe is the sequence 
$\as{0}...$ of addresses of its iterates, defined by
\begin{equation}\label{eqn:address}
	\ftoof{i}{x}\in\Shiftset{\as{i}}.
\end{equation}
However, when $\kerH\neq\emptyset$,  \refer{eqn}{address} need not define a unique itinerary for a point \inX{x}.  Let \cylinders{n} denote the set of $n$-tuples (or ``$n$-words'') $\as{0}...\as{n-1}$ with $\as{i}\in\single{1,...,m}$.
For \inNat{n}, the \deffont{$\boldsymbol{n}$-itinerary set} of \inX{x} is the subset of \cylinders{n}
defined via \refer{eqn}{address}:
\begin{equation*}
	\nitinset{x}{n}\eqdef
		\setbld{\ads{0}\dots\ads{n-1}\in\cylinders{n}}
			{\ftoof{i}{x}\in\Shiftset{\ads{i}}\text{ for }i<n}
		={\bigtimes}_{i=0}^{n-1}\adset{\ftoof{i}{x}}.
\end{equation*}
The set of points for which $\nitinset{x}{n}\neq\emptyset$
\begin{equation*}
	\DomsH{n}\eqdef\bigcap_{i=0}^{n-1}\fpre{i}{\Shiftpiece}
\end{equation*}
is the union of the sets
\begin{equation}\label{eqn:pi}
	\realizeset{w}\eqdef\setbld{\inX{x}}{w\in\nitinset{x}{n}}
		=\bigcap_{i=0}^{n-1}\fpre{i}{\Shiftset{\ws{i}}}
\end{equation}
as $w=\ws{0}...\ws{n-1}$ ranges over the finite collection \cylinders{n} of $n$-words. The intersection
	$$\DomH\eqdef\bigcap_{i=0}^{\infty}\fpre{i}{\Shiftpiece}=\bigcap_{n=1}^{\infty}\DomsH{n}$$
is the set of points whose whole forward orbit is contained in \Shiftpiece{}; we refer to it as the \deffont{domain} of \Shiftsystem.

For a horseshoe, the assignment of an itinerary to each point of \DomH{} is a continuous map onto the $m$-shift space \mshiftspace{} and hence provides a semiconjugacy from the \DomH{}-restriction of $f$ to the one-sided shift map on \mshiftspace{};  the estimate \eqref{eqn:main} on entropy is an immediate consequence in this case.

The continuity of itineraries with respect to a horseshoe has a semicontinuity analogue for a general closed shift system.  Recall that for any sequence of sets \As{i}, \inflist{i}{1,2}
\begin{equation*}
	\limsup\As{i}\eqdef\bigcap_{k=1}^{\infty}\bigcup_{i=k}^{\infty}\As{i}.
\end{equation*}

\begin{lemma}\label{lem:semicont}
	Suppose \Shiftsystem{} is a closed $m$-shift system.
	\begin{enumerate}
		\item \label{eqn:semicont1} For each \inX{x}, \itinset{x} is a closed subset of \mshiftspace{}.
		
		\item \label{eqn:semicont2} For each nonempty (closed) set $A\subset\cylinders{n}$,
			\inNat{n} ($n=\infty$), \realizeset{A} is a nonempty closed subset of \X{}.
		
		\item \label{eqn:semicont3} For \inNat{n}, if \inDomsH{n}{\xs{i}} for \inflist{i}{1,2}, then
			$\limsup\nitinset{\xs{i}}{n}\neq\emptyset$.
		
		\item \label{eqn:semicont4} The set-valued maps $x\mapsto\nitinset{x}{n}$,
			$n\in\Nat\cup\single{\infty}$, are upper semicontinuous: if $\xs{i}\to x$
			in \X{}, then $\limsup\nitinset{\xs{i}}{n}\subset\nitinset{x}{n}$.
	\end{enumerate}
\end{lemma}
This is Lemma 3.3 of \cite{BobokZN:mfold};  we refer the reader there for a proof.

To obtain entropy estimates from a shift system with nonempty kernel, we need some further assumptions.

\begin{definition}\label{dfn:locdiv}
	We say a set $W\subset\X$ \deffont{locally divides} the $m$-shift system{} \Shiftsystem{} if
	\begin{enumerate}
		\item\label{locdiv1} $\fof{W}\subset W$
		\item\label{locdiv2} $\Shiftset{i}\setminus W\nonempty$ for \ilist{i}{1}{m}
		\item\label{locdiv3} There exists a closed shift-invariant set $\Lambda\subset\mshiftspace$
			such that 
				$$\ent{\Lambda}<\log m$$
			and a neighborhood $V$ of $W$ in \X{} such that if $\ftoof{i}{x}\in V\setminus W$ 
			for \ilist{i}{0}{n-1} then 
				$$\nitinset{x}{n}\subset\Lambda(n)$$
			where $\Lam(n)$ denotes the set of initial words of length $n$ for sequences in \Lam{}.
	\end{enumerate}
\end{definition} 
We note that the third condition above is somewhat different from that in the definition of local division in \cite{BobokZN:mfold}, but contains it as a special case by \cite[Lemma 4.2]{BobokZN:mfold}.  


\begin{markedremark}\label{rmk:nonempty}
	Let \Shiftsystem{} be a closed $m$-shift system, and $w\in\cylinders{n}$ any finite word.  Then
	\begin{enumerate}
		\item \label{nonempty1} 
			\realizeset{w}  is a closed nonempty set in \X{}
		\item \label{nonempty2} If $W$ locally divides \Shiftsystem, then  \realizeset{w}
			is not contained in $W$.
	\end{enumerate}
\switchtotext
	\eqref{nonempty1} is an easy consequence of \refer{eqn}{pi}, while
	\eqref{nonempty2} is a consequence of conditions \ref{locdiv1} and \ref{locdiv2} in 
	\refer{dfn}{locdiv}.
\end{markedremark}

\begin{lemma}\label{lem:unionlocdiv}
	If \Ws{j}, \ilist{j}{1}{n} are sets that locally divide a closed $m$-shift system \Shiftsystem{}, 
	then their union 
	$W\eqdef\bigcup_{j=1}^{n}\Ws{j}$ also locally divides 
	\Shiftsystem{}, provided that $\Shiftset{i}\setminus W\nonempty$ for \ilist{i}{1}{m}.
\end{lemma}

\begin{proof}
	The only condition from Definition \ref{dfn:locdiv} which is not immediate is (3).
	Let \Lams{j} and \Vs{j} be the shift-invariant set and neighborhood specified 
	in \subref{dfn}{locdiv}{locdiv3} for \Ws{j}, and set
		$$\Lamp\eqdef\bigcup_{j=1}^{n}\Lams{j}.$$
	Clearly $\ent{\Lamp}=\max\ent{\Lams{j}}<\log m$, and we can choose 
	a word $w=\ws{0}...\ws{k-1}$ of (some) length $k$ 
	which does not appear in any sequence in \Lamp{}:
		$$w\in\mshiftspace(k)\setminus\Lamp(k).$$
	Pick $i\neq\ws{k-1}\in\setlist{1}{m}$.
	
	Since each \Ws{j} is invariant, we can find a neighborhood \Gs{j} of \Ws{j} with 
		$$\bigcup_{i=0}^{k-1}\ftoof{i}{\Gs{j}}\subset\Vs{j}.$$
	Then 
		$$V=\bigcup_{j=1}^{n}\Gs{j}$$
	is a neighborhood of $W$. For \inflist{\ell}{1}, set
		$$\Fs{\ell}\eqdef\bigcap_{i=0}^{\ell-1}\fpre{i}{V\setminus W}$$
	and define $\Ftil{\ell}\subset\mshiftspace$ to consist of all words of the form $\alpha iii....$, where
	$\alpha$ belongs to the \el{}-itinerary set of some point in \Fs{\ell}.
	
	Now the set
		$$\Lam\eqdef \single{iii...}\cup\clos\bigcup_{\ell\geq1}\Ftil{\ell}$$
	is clearly a closed shift-invariant subset of \mshiftspace{}.  Moreover, if $\ftoof{j}{x}\in V\setminus W$
	for \ilist{i}{0}{\ell-1}, then $\nitinset{x}{\ell}\subset\Lam(\ell)$.  To complete the proof of the lemma
	it suffices to show $w\notin\Lam(k)$, which implies that $\ent{\Lam}<\log m$.
	
	If $w\in\Lam(k)$, then it belongs to $\nitinset{\Fs{\ell}}{k}$ for some $\ell$.  But then 
	by the choice of $V$, $w$ belongs to \nitinset{x}{k} for some point $x$ with 
	$\ftoof{i}{x}\in V_{j}\setminus W_{j}$, \ilist{i}{0}{k}, for some $j$.  Hence 
	$w\in\Lams{j}(k)\subset\Lamp(k)$, contrary to our choice of $w$.
\end{proof}

The following is a modification of \cite[Lemma 4.4]{BobokZN:mfold} to fit our more general definition of local division.

\begin{lemma}\label{lem:locdivmin}
	Suppose \Shiftsystem{} is a closed $m$-shift system 
	and $W$ is a set which locally divides \Shiftsystem{} and contains all minimal sets
	in the core \centOH{} (\ie{} which are contained in the center \centH{}).
	
	Then there exists \inNat{\zeta} such that any orbit segment of length $\zeta$ which is contained
	in \centH{} terminates in $W$.
\end{lemma}

\begin{proof}
	Suppose $\{\ftoof{j}{\xs{n}}\}_{j=0}^{n-1}$ are (arbitrarily long) orbit segments contained in 
	$\centH\setminus W$ and (passing to a subsequence if necessary) assume $\xs{n}\to x$.
	Then the orbit of $x$ is contained in \centH{} and so there is a minimal set
	$M\subset\omega(x)\cap W$. In particular, the continuity of $f$ implies that we can find
	orbit segments $\{\ftoof{j}{x}\}_{j=k_{n}}^{k_{n}+n-1}$ of $x$ with increasing length 
	contained in $V$, and hence we can find a subsequence $\{\xs{l_{n}}\}$ of our original points
	whose orbit segments $\{\ftoof{j}{\xs{l_{n}}}\}_{j=k_{n}}^{k_{n}+n-1}$ are contained in
	$\centH\cap(V\setminus W)$.  But this means that 
	$\nitinset{\ftoof{k_{n}}{x_{l_{n}}}}{n}\subset\Lam(n)$ consists of
	all words of length $n$, contradicting condition (3) of Definition \ref{dfn:locdiv}.
\end{proof}

Given a closed $m$-shift system \Shiftsystem{} for \fXX{}, we can  associate to any closed $f$-invariant set $S\subset\DomH$ two different ``entropies'':  the topological entropy of the restriction of $f$ to $S$, $\entof{S}\eqdef\ent{f|S}$, and the topological entropy of the restriction of the shift map to the itinerary set of $S$, \vent{S}, which we refer to as the \deffont{virtual entropy} of $S$.  These are not related in any \apriori{} way--in particular the virtual entropy of $f$ on a periodic orbit, unlike the topological entropy, need not be zero.
However, we can sometimes get an \apriori{} bound on it.

\begin{prop}\label{prop:locdivide}
	Under the conditions of \refer{lem}{locdivmin}, there exists $\beta<\log m$ such that any periodic
	orbit $P\subset\DomH$ which is not contained in the set $W$ (and hence is disjoint from the core)
	has virtual entropy bounded by $\beta$:
		$$\vent{P}\leq\beta.$$
\end{prop}
The proof of this is the same as that of \cite[Proposition 4.5]{BobokZN:mfold}, with the center \centH{} replaced by the set $W$.


\begin{lemma}\label{lem:mincount}
	Suppose \Shiftsystem{} is a closed $m$-shift system 
	whose kernel is \deffont{eventually countable}:
		$$\ftoof{j}{\kerH}\text{ is (at most) countable for some }j.$$
		
	Then for any infinite minimal set $M\subset\X$,
		$$\vent{M}\leq\ent{M}.$$
\end{lemma}
The proof of this is the same as that of \cite[Lemma 4.7]{BobokZN:mfold}.

From this, we have the following analogue of \cite[Theorem 4.8]{BobokZN:mfold}:

\begin{theorem}\label{thm:mshiftent}
	Suppose \fXX{} has a closed, nontrivial $m$-shift system \Shiftsystem{} for which
		\begin{enumerate}
			\item the kernel is eventually countable 
			\item there exists a set $W\subset\X$ such that
			\begin{enumerate}
				\item $W$ contains all minimal sets in the core \centOH{}, and 
				\item $W$ locally divides \Shiftsystem{}.
			\end{enumerate}
		\end{enumerate}
		
		Then
			$$\ent{f}\geq\log m.$$
\end{theorem}
Note that if the kernel is eventually countable, then every minimal set in the core \centOH{} is a periodic orbit.  

To prove Theorem \ref{thm:mshiftent}, we first establish some preliminary results.

\begin{remark}\label{rmk:itin}
	Suppose $x\in X$ and $A\subset \mshiftspace$ satisfy, for some $j,k\in\mathbb{N}$
	\begin{equation*}
		\nitinset{\ftoof{j}{x}}{k}\cap\sigma^{j}[A](k)=\emptyset.
	\end{equation*}
	Then
	\begin{equation*}
		\nitinset{x}{k+j}\cap A(k+j)=\emptyset.
	\end{equation*}
\end{remark}
The proof of this is the same as that of \cite[Remark 4.9]{BobokZN:mfold}.

\begin{lemma}\label{lem:mshiftdisj}
	Suppose $A\subset\mshiftspace$ is closed and shift-invariant, and $y\in\omega(x)\subset X$
	with $\Omega(y)\cap A=\emptyset$.
	
	Then there exists a neighborhood $U$ of $x$ and $k\in\mathbb{N}$ such that 
	\begin{equation*}
		\nitinset{x^{\prime}}{k}\cap A(k)=\emptyset
	\end{equation*}
	for every $x^{\prime}\in U$.
\end{lemma}
This is the same as \cite[Lemma 4.10]{BobokZN:mfold}.

\begin{prop}\label{prop:msep}
	Suppose \Shiftsystem{} satisfies the hypotheses of Theorem \ref{thm:mshiftent}, and that
	$\Gamma\subset\mshiftspace$ is a shift-minimal set such that, for every $f$-minimal set
	$M$ disjont from $W$,
	\begin{equation*}
		\entof{\Gamma}>\max\{\entof{\Lam}. \entof{\Omega(M)}\}
	\end{equation*}
	where \Lam{} is as in Definition \ref{dfn:locdiv}.
	
	Then there exists $k\in\mathbb{N}$ such that every point $x\in X$ with 
	$\nitinset{x}{k}\cap\Gamma(k)\neq\emptyset$ satisfies $\ftoof{k-1}{x}\in W$.
\end{prop}

\begin{proof}
	We construct for every $x\in X$ a neighborhood $U(x)$ and an associated integer $k(x)$ such
	that every point $x^{\prime}\in W$ with $\ftoof{k(x)-1}{x^{\prime}}\notin W$ has
	$\nitinset{x^{\prime}}{k(x)}\cap\Gamma(k(x))=\emptyset$.  We consider three cases;  even though
	the second and third need not be mutually exclusive, this presents no problem:
	\begin{enumerate}
		\item If $x\notin\DomH$, pick $k(x)$ so that $\ftoof{k(x)}{x}\notin\Shiftpiece$, 
		and a neighborhood $U(x)$ of $x$ for which $\ftoof{k(x)}{U(x)}\cap\Shiftpiece=\emptyset$.
		Then $\nitinset{x^{\prime}}{k(x)}=\emptyset$ for all $x^{\prime}\in U(x)$.
		
		\item If $\omega(x)$ contains a minimal set $M$ which is not contained in $W$, then since 
		$\entof{\Omega(M)}<\entof{\Gamma}$ and \Gam{} is minimal, $\Omega(M)$ is disjoint
		from \Gam{}.  It follows by Lemma \ref{lem:mshiftdisj} 
		with $y$ any element of $M$ and $A=\Gamma$ that we can find $U(x)$ and $k(x)$ so that
		$\nitinset{x^{\prime}}{k(x)}\cap\Gam(k(x))=\emptyset$ for all $x^{\prime}\in U(x)$.
		
		\item If $M\subset\omega(x)\cap W$, pick \Lam{} and $V$ associated to $W$ in Definition
		\ref{dfn:locdiv}.  Since $\entof{\Lam}<\entof{\Gam}$, \Lam{} is disjoint from \Gam{}, and 
		hence $\Lam(k_{0})\cap\Gam(k_{0})=\emptyset$ for some $k_{0}\in\mathbb{N}$.
		But since $M\subset\omega(x)$, there exists $k_{1}\in\mathbb{N}$ so that
		$\ftoof{k_{1}+j}{x}\in V$ for $0\leq j<k_{0}$, and a neighborhood $U(x)$ so that the
		same holds true for every $x^{\prime}\in U(x)$.  Let $k(x)=k_{0}+k_{1}$.  For any 
		$x^{\prime}\in U(x)$ with $\ftoof{k(x)-1}{x^{\prime}}\notin W$ we have
		$\ftoof{k_{1}+j}{x^{\prime}}\in V\setminus W$ for $0\leq j<k_{0}$, and hence
		\nitinset{x^{\prime}}{k(x)} is disjoint from $\Gam(k(x))$ by Lemma \ref{lem:mshiftdisj}.
	\end{enumerate}
	
	Since $\omega(x)$ always contains some minimal set, these cases are exhaustive, and so
	\setbld{U(x)}{x\in X} form an open cover of $X$.  leet \setbld{U(\xs{i}}{\ilist{i}{1}{N}} be a 
	finite subcover, and set
	\begin{equation*}
		k=\max_{\ilist{i}{1}{N}}k(\xs{i}).
	\end{equation*}
	Then we clearly have the desired conclusion with this value of $k$.
\end{proof}

\begin{proof}[Proof of \refer{thm}{mshiftent}]
%
	Let \betha{} as in \refer{prop}{locdivide}. 
	We will show that for $0<\eps<\log m-\max\single{\entof{\Lam},\beta}$, $f$ has minimal sets $M$
	with 
	\begin{equation*}
		\ent{M}\geq\log m-\eps.
	\end{equation*}
	By \cite{Grillenberger}, \mshiftspace{} contains shift-minimal sets with entropy arbitrarily near
	$\log m$, so we can find \Gams{\eps} minimal with 
	\begin{equation*}
		\entof{\Gams{\eps}}>\log m-\eps>\max\single{\entof{\Lam},\beta}.
	\end{equation*}
	If some $M\subset\Shiftpiece\setminus W$ has $\ent{M}<\log m-\eps$, 
	then by \refer{prop}{locdivide}
	and \refer{lem}{mincount}, $\entof{\Omega(M)}<\log m-\eps$.  Thus, if no minimal set 
	$M$ has $\entof{\Omega(M)}\geq\log m-\eps$ , then \refer{prop}{msep} says that for some 
	\inNat{k}, $\nitinset{x}{k}\cap\Gams{\eps}(k)=\emptyset$ whenever $\ftoof{k-1}{x}\notin W$.
	But \refer{rmk}{nonempty}(3) says that every $w\in\Gams{\eps}(k)$ belongs to some
	\nitinset{x}{k} for a point with $\ftoof{k-1}{x}\notin W$, a contradiction.
	
	This establishes the existence of minimal sets satisfying $\ent{M}\geq\log m-\eps$.  Thus
	$\ent{f}\geq\log m-\eps$, and since \epsgo{} can be chosen arbitrarily small, the conclusion follows.
\end{proof}


\section{Trees}\label{sec:background}

Topologically, a (finite) \deffont{tree} is a uniquely arcwise connected Hausdorff space which is a union of (finitely many) closed intervals.  The complement $\T\setminus\single{x}$ of a point $x\in\T$ has finitely many components, called the \deffont{branches} of \T{} at $x$;  a \deffont{closed branch} at $x$ is the union of \single{x} with a branch at $x$.  The \deffont{valence} of $x$ is the number of branches at $x$;  a point of valence one \resp{valence $>2$} is called an \deffont{endpoint} \resp{a \deffont{branchpoint}} of \T{};  the set of all branchpoints of \T{} is denoted \Branchpts{}.  

We endow \T{} with further combinatorial structure, first by distinguishing a finite set \Vertices{} of \deffont{vertices} which includes all endpoints and branchpoints (and perhaps some valence two points), and then distinguishing one vertex \vs{0} as the \deffont{root} of \T{}.  An \deffont{edge} of \T{} is the closure of a component of $\T\setminus\Vertices$.  Each closed branch of \T{} at $v\in\Vertices$ contains a unique edge with endpoint $v$; for $v\neq\vs{0}$, the \deffont{incoming} branch \resp{incoming edge} at $v$ is the 
branch containing \vs{0} 
\resp{the unique edge at $v$ contained in the convex hull \hull{\vs{0},v}}, and the other branches at $v$--as well as edges at $v$ contained in their closures--are \deffont{outgoing} at $v$.  We direct each edge of \T{} so that $v$ is the terminal or right \resp{initial or left} endpoint of the incoming \resp{any outgoing} edge at $v$.  We use interval notation, denoting the edge with left endpoint $u$ and right endpoint $v$ by \clint{u}{v}, and adapt the notation of open and half-open intervals to denote edges missing one or both endpoints.  
The number of outgoing branches (equivalently edges) at $v\in\Vertices$ is its \deffont{outdegree} (clearly, for $v\neq\vs{0}$ this is one less than the valence). The \deffont{level} of a vertex $v\in\Vertices$ is the number of edges contained in \hull{\vs{0},v} (so the root is at level zero).  

Now, we wish to define a linear ordering \admitless{} on the points of \T{}.  We begin by numbering the 
outgoing edges at each vertex $v\in\Vertices$; this induces a numbering of the outgoing branches,  \Branchat{i}{v}, \ilist{i}{1}{\valence}, where \valence{} is the outdegree of $v$;  for the moment, this numbering is arbitrary, but we will
impose a further condition on it in \refer{sec}{center}.  For $v\neq\vs{0}$, the incoming branch is numbered zero, so 
	$$\T\setminus\single{v}=\bigcup_{i=0}^{\valence}\Branchat{i}{v}$$
(for $v=\vs{0}$, the only difference is that there is no \Branchat{0}{\vs{0}}).  

Given this numbering, we assign to each point $x\neq\vs{0}$ an address as follows.  There is a unique simple path \gamof{x} from \vs{0} to $x$; let $\Vof{x}=(\vs{0},...,\vs{k})$ be the sequence of vertices occurring along \gamof{x}---if $x\in\Vertices$ then \vs{k} is the last vertex along \gamof{x} before $x$.  
For each \ilist{j}{0}{k}, there is precisely one outgoing branch at \vs{j}, say \Branchat{\is{j}}{\vs{j}}, containing $x$.  The sequence $\boldsymbol{\alpha(x)}\eqdef(\is{0},...,\is{k})$ is the \deffont{address} of $x$ in \T{}. Two points have the same address in \T{} precisely if they belong to the same left-open edge \lopint{\vs{k}}{\vs{k+1}}.  The linear ordering \admitless{} is then defined by lexicographic comparison of addresses and, within an edge, the direction from left to right.  More precisely:

\begin{definition}[The Linear Ordering \admitless{} on \T{}]\label{dfn:admitless}
	\begin{enumerate}
		\item $\vs{0}\admitless x$ for every $x\neq\vs{0}$ in \T{}.
		\item Given $x\neq\xp$ in \T{} with $\alpha(x)=(\is{0},...,\is{k})$ 
		and $\alpha(\xp)=(\ips{0},...,\ips{\kp})$,
			\begin{enumerate}
				\item if $\alpha(x)=\alpha(\xp)$, then $x$ and \xp{} both belong to a
					common left-open edge \lopint{\vs{k}}{\vs{k+1}};  we write $x\admitless\xp$
					if $x\in\hull{\vs{k},\xp}$;
				\item if $\is{j}=\ips{j}$ for \ilist{j}{0}{k} (and $\kp>k$) then $x\admitless\xp$;
				\item if $\js{0}\eqdef\min\setbld{j}{\is{j}\neq\ips{j}}$ then 
					$x\admitless\xp$ iff $\is{\js{0}}<\ips{\js{0}}$.
			\end{enumerate}

	\end{enumerate}

\end{definition}

For $v\in\Vertices$ with outdegree \valence{} and any $i\in\single{1,...,\valence}$, set
\begin{align*}
	\Downbranches{i}{v}&\eqdef\bigcup_{0<\ip<i}\Branchat{\ip}{v}\\
	\Upbranches{i}{v}&\eqdef\bigcup_{\ip>i}\Branchat{\ip}{v}.
\end{align*}

\begin{remark}[Topological interpretation of the ordering \admitless{} on \T{}]\label{rmk:interp}
 	\begin{enumerate}
		\item If $\xp\in\interior\hull{\vs{0},x}$ then $\vs{0}\admitless\xp\admitless x$.

		\item For any vertex $v\in\Vertices$ with $\Vof{v}=(\vs{0},...,\vs{k})$ and 
			$\alpha(v)=(\is{0},...,\is{k})$,
			\begin{align*}
				\setbld{x}{x\admitmore v}
					&=\bigcup_{j=0}^{k}\Upbranches{\is{j}}{\vs{j}}\cup\bigcup_{i>0}\Branchat{i}{v}\\
				\setbld{x}{x\admitless v}
					&=\bigcup_{j=0}^{k}\Downbranches{\is{j}}{\vs{j}}
						\cup\hull{\vs{0},v}\setminus\single{v}.
			\end{align*}
	
		\item For $x\in\opint{\vs{k}}{\vs{k+1}}$ with $\Vof{x}=(\vs{0},...,\vs{k})$ 
			and $\alpha(x)=(\is{0},...,\is{k})$,
			\begin{align*}
				\setbld{\xp}{\xp\admitmore x}
					&=\lopint{x}{\vs{k+1}}\cup\setbld{\xp}{\xp\admitmore\vs{k+1}}\\
					&=\lopint{x}{\vs{k+1}}\cup
						\bigcup_{j=0}^{k}\Upbranches{\is{j}}{\vs{j}}
						\cup\bigcup_{i>0}\Branchat{i}{\vs{k+1}}\\
				\setbld{\xp}{\xp\admitless x}
					&=\ropint{\vs{k}}{x}
						\cup\Downbranches{\is{k}}{\vs{k}}
						\cup\setbld{\xp}{\xp\admitless\vs{k}}\\
					&=\bigcup_{j=0}^{k}\Downbranches{\is{j}}{\vs{j}}
						\cup\hull{\vs{0},x}\setminus\single{x}.
			\end{align*}
	\end{enumerate}
\end{remark}

	Note in particular that for any $v\in\Vertices{}$ 
	the \emph{outgoing} branches at $v$ are comparable:
	for $0<i<\ip$ if $x\in\Branchat{i}{v}$ and $\xp\in\Branchat{\ip}{v}$ then 
	$x\admitless\xp$;  however this is in general false if $i=0$.
	
	A key property of the linear ordering on the real line is its continuity:  that if two 
	convergent sequences $\xs{i}\to x$ and $\xps{i}\to\xp$
	satisfy $\xs{i}\leq \xps{i}$ for all $i$, then their limits satisfy the same inequality:
	$x\leq\xp$.  This is false for our ordering: if \vs{j} is a vertex belonging to \Vof{\xp}
	and \is{j} is the corresponding element of $\alpha(x)$, then a sequence 
	$\xps{i}\in\Upbranches{\is{j}}{\vs{j}}$ will satisfy $x\admitless\xps{i}$, but if
	it converges to $\xp=\vs{j}$ then $\lim\xps{i}\admitless x$.  However, this can 
	only happen when the limit is a vertex.

\begin{lemma}\label{lem:orderlim}
	Suppose $\xs{i}\to x$ and $\xps{i}\to\xp$ are convergent sequences in \T{} with
	$$\xs{i}\admitless\xps{i},\quad\inflist{i}{1,2}.$$
	If $x,\xp\notin\Vertices$, then either $x=\xp$ or 
	$$x\admitless\xp.$$
\end{lemma}

\begin{proof}
	If $x$ and \xp{} are interior to the same edge of \T{}, the conclusion is trivial.  So 
	suppose not.  Then they are interior to distinct edges of \T{}, and these are comparable.
	Since the convergent sequences are eventually interior to the corresponding edges, the
	conclusion is immediate.
\end{proof}

\section{The Kernel}\label{sec:kernel}
In this section, we show, given a \comfold{} map \fTT{} on the tree \T{}, how to construct an 
$m$-shift system \Shiftsystem{} with eventually countable kernel (in fact, with \fof{\kerH} countable), thus fulfilling the first condition of \refer{thm}{mshiftent}.  This is based on the idea of \emph{regular $m$-sections} from \cite{BobokZN:mfold}.
 
 

\subsection{Regular values}\label{subsec:reg}
Suppose \fTT{} is continuous and \inT{x} with $y=\fof{x}\in\T\setminus\Vertices$.  Then a small neighborhood of $y$ is an interval disjoint from \Branchpts{}, so that it makes sense to talk about points being on one or the other side of $y$.  We say that $x$ is a \deffont{non-minimal} \resp{\deffont{non-maximal}} preimage of $y$ if there exist points \xp{} arbitrarily near $x$ with $\fof{\xp}<y$ \resp{$\fof{\xp}>y$}.



Given \inNat{m}, we say that $y\in\T\setminus\Vertices$ is a \deffont{left $\boldsymbol{m}$-regular} \resp{\deffont{right $\boldsymbol{m}$-regular}} value for \fTT{} if $y$ has at least $m$ non-minimal \resp{non-maximal} preimages.  The set of all \leftreg{m} \resp{\rightreg{m}} values of $f$ will be denoted \Lregm{} \resp{\Rregm}, and their intersection, the set of \deffont{$\boldsymbol{m}$-regular values}, will be denoted \Regm{}.

In \cite[\S5.1]{BobokZN:mfold} it is shown that \Regm{} is an open cocountable set for any \comfold{} map of the circle;  the analogous result in our case, \fTT{}, \T{} a (finite) tree, follows from an analogous argument.

\begin{prop}\label{prop:regopen}
	For any \comfold{} map \fTT{} on a finite tree, \Regm{} is an open cocountable set.
\end{prop}

\begin{proof}
	The argument in \cite{BobokZN:mfold} has four steps;  we summarize these and indicate any 
	adjustments to make each step work in our setting:
	\begin{enumerate}
		\item \cite[Lemma 5.1]{BobokZN:mfold}  \textit{If \Interval{} and $\Y=\fof{\Interval}$ are 
		closed intervals,
		then every $y\in\interior\Y$ has at least one non-minimal 
		and at least one non-maximal preimage
		in \Interval{}.
		}
		
		This is because the components of $\Interval\setminus\fprey$ containing the endpoints 
		of \Interval{} must map onto one-sided neighborhoods of $y$.
		
		\item \cite[Lemma 5.2]{BobokZN:mfold}  \textit{If $y$ is left 
		\resp{right} \reg{m}, then the 
		interior of some left \resp{right} neighborhood of $y$ is contained in \Regm{}.
		}
		
		This carries over, since we can separate $m$ non-minimal \resp{non-maximal} preimages
		of $y$ with neighborhoods disjoint from \Branchpts{}.
		
		\item \cite[Lemma 5.3]{BobokZN:mfold}  \textit{If \fprey{} has at least $m$ components, then
		there exists a nontrivial interval \Y{} with $y$ an endpoint such that $\interior\Y\subset\Regm$.
		}
		
		This carries over provided $y\notin\fof{\Vertices}$, since we can find $2m$ intervals, each 
		contained in a single edge of \T{}, with one endpoint an endpoint of a component of \fprey{},
		and each mapping onto a one-sided neighborhood of $y$;  at least $m$ of these intervals map
		onto the same side of $y$, and we can apply the preceding result.
		
		\item \cite[Proposition 5.4]{BobokZN:mfold}  \Regm{} is open by the second result above,
		and the points $y$ for which \fprey{} has nonempty interior is at most countable.  
		Throwing these out as well as the (finite) set \fof{\Vertices}, we have a \cocount{} set 
		$\Y\subset\T$ for which the third result above says that each \inY{y} is either an element
		of \Regm{} or an endpoint of a component of \Regm{}.  Since \Regm{} is open, it has
		(at most) countably many components, so throwing away their endpoints (from \Y{}) we
		obtain a cocountable subset of the open set \Regm{}.
	\end{enumerate}

\end{proof}

\subsection{Regular sections}\label{subsec:regsect}

If \fTT{} is \mfold{} on the subset \Y{}, then an  \deffont{$\boldsymbol{m}$-section} for $f$ on \Y{} is a choice for each \inY{y} of a collection of $m$ distinct preimages $\psijof{y}\in\fprey$, \ilist{j}{1}{m}, which we combine as a map
	$$\psi:\Y\to\bigtimes_{i=1}^{m}\T\eqdef\T\times...\times\T.$$

We say $\psi$ is a \deffont{left-regular section} \resp{\deffont{right-regular section}} if \psisof{i}{y} is a
non-minimal \resp{non-maximal} preimage of $y$ for every \inY{y}.
Clearly, a left-regular \resp{right-regular} section can only be defined on a subset of \Lregm{} \resp{\Rregm}, but even if $\Y\subset\Regm$ it need not be possible to define a section on \Y{} which is both left-regular and right-regular.

We define several calibrations of the ``spread'' of an \msec{} $\psi$ on $\Y\subset\T$.
Given \inY{y}, set
	$$\deltof{\psi}{y}\eqdef\min_{1\leq j<j\pr\leq m}\distof{\psijof{y}}{\psisof{j\pr}{y}}$$
and then the \emph{mesh} of $\psi$ on the set $U\subset\Y$ is
	$$\Delof{\psi}{U}\eqdef\inf\setbld{\deltof{\psi}{y}}{y\in U}.$$
Also, the distance on the product $\bigtimes_{{i=1}}^{m}\T$ defines a distance between the values of $\psi$ at two points \inY{y,y\pr}
	$$\mdistof{\psiof{y}}{\psiof{y\pr}}\eqdef\max_{\ilist{j}{1}{m}}\distof{\psijof{y}}{\psijof{y\pr}}$$
and the \emph{variation} of $\psi$ across a set $U\subset\Y$ is
	$$\varof{\psi}{U}\eqdef\sup\setbld{\mdistof{\psiof{y}}{\psiof{y\pr}}}{y,y\pr\in U}.$$

We refer the reader to \cite{BobokZN:mfold} for the proof of the following result.  

\begin{prop}[Proposition 5.8, \cite{BobokZN:mfold}]\label{prop:regsec}
	Suppose \fTT{} is a cocountably \mfold{} map. 
	Then there is a \cocount{} open subset $Y\subset\Regm$ and a left-regular
	\resp{right-regular} \msec{} $\psi$ of $f$ on \Y{} such that, for every component \Ys{i}
	of \Y{}
	\begin{equation}\label{eqn:regsec1}
		\Delof{\psi}{\Ys{i}}>0
	\end{equation}
	\begin{equation}\label{eqn:regsec2}
		\varof{\psi}{\Ys{i}}<\half\Delof{\psi}{\Ys{i}}.
	\end{equation}
\end{prop}



%
%
The condition \eqref{eqn:regsec2} has several useful consequences.

\begin{lemma}\label{lem:shift}
	Suppose that, as in \refer{prop}{regsec}, $\psi$ is an \msec{} of $f$ on the 
	\cocount{} open set $\Y\subset\T$ such that every component \Ys{i} of \Y{} satisfies
	\refer{eqn}{regsec2}.
	
	If for some sequence \inY{\ys{i}} and some sequence \inZ{\js{i}} we have
	\psisof{\js{i}}{\ys{i}} converging to a point $x\in\fpreof{\Y}$, then \js{i} is eventually
	constant. 
\end{lemma}


\begin{proof}
	Suppose $\js{1}\neq\js{2}$ and 
		\inY{\ys{i},\yps{i}} satisfy
			$$u=\lim\psisof{\js{1}}{\ys{i}}=\lim\psisof{\js{2}}{\yps{i}}$$
		with \inY{\fof{u}}.
		By continuity of $f$,
			$$y\eqdef\fof{u}=\lim\fof{\psisof{\js{1}}{\ys{i}}}=\lim\ys{i}=\lim\yps{i}.$$
		We can assume that $y$, \ys{i} and \yps{i} all belong to the same component $U$ of \Y{}.
		For $i$ large,
			$$\distof{\psisof{\js{1}}{\ys{i}}}{x}<\frac{1}{4}\Delof{\psi}{U}\text{ and }
			\distof{\psisof{\js{2}}{\yps{i}}}{x}<\frac{1}{4}\Delof{\psi}{U}.$$
		But then
			\begin{eqnarray*}
				\Delof{\psi}{U} &\leq \distof{\psisof{\js{1}}{\ys{i}}}{ \psisof{\js{2}}{\ys{i}}}\\
					&\leq \distof{\psisof{\js{1}}{\ys{i}}}{ \psisof{\js{2}}{\yps{i}}}
						+\distof{\psisof{\js{2}}{\yps{i}}}{\psisof{\js{2}}{\ys{i}}}\\
					&\leq\distof{\psisof{\js{1}}{\ys{i}}}{x}
						+\distof{x}{\psisof{\js{2}}{\yps{i}}}
						+\distof{\psisof{\js{2}}{\yps{i}}}{\psisof{\js{2}}{\ys{i}}}\\
					&\leq\varof{\psi}{U}+\half\Delof{\psi}{U}\\
					&<\Delof{\psi}{U},
			\end{eqnarray*}
		a contradiction.
\end{proof}

\begin{remark}\label{rmk:renumber}
	If an \msec{} $\psi$ satisfies \refer{eqn}{regsec2}, and we modify $\psi$ by applying a permutation
	to the indices \ilist{j}{1}{m} of the components \psisof{j}{y} for all $y$ in some component \Ys{i} of
	\Y{}, then \refer{eqn}{regsec2} remains true. 
\end{remark}


\subsection{Monotone Sections}\label{sec:monsec}

 Suppose that we have a linear ordering $\admitless$ on \T{} as in \refer{sec}{background}.  
We begin with some remarks on comparability.

\begin{markedremark}[Comparability]\label{rmk:comparable}
	\begin{enumerate}
		\item Suppose $x\admitless\xp$ and $x$ is not a vertex.  Then there exist neighborhoods
		$U,U\pr$ of $x$ and $\xp$ which
		are comparable:  that is, $z\admitless \zp$ for all $z\in U$ and $\zp\in U\pr$,
		and we write $U\admitless U\pr$.
		
		\item An interval disjoint from \Branchpts{} is comparable with any point not contained in it.
	\end{enumerate}
	\switchtotext
	This is an immediate consequence of \refer{rmk}{interp}.
\end{markedremark}

\begin{markedremark}\label{rmk:intervals}
		Suppose $J$ is a nontrivial open interval contained in $\fof{\hull{z,\zp}}$ 
		and disjoint from $\Branchpts$.
		Then there exists an open interval $J\pr\subset\hull{z,\zp}$ with $\fof{J\pr}=J$.
	\switchtotext
		To see this, consider the preimage sets for the two endpoints of $J$; these are
		disjoint nonempty closed subsets of the interval \hull{z,\zp}, and it is easy to see that in such 
		a situation some component of the complement of their union is an interval with one endpoint 
		in each, and it must map onto $J$.
\end{markedremark}

We say that $\psi$ is a \deffont{monotone section} if $\psisof{i}{y}\admitless\psisof{i+1}{y}$ for each \ilist{i}{1}{m-1} and $y\in\Y$.  

\begin{markedremark}\label{rmk:monsec}
		Suppose $f$ has an \msec{} $\psi$ defined on a cocountable open set $\Y\subset\T$
		satisfying \refer{eqn}{regsec2}.  Then we can renumber the indices of the components
		by a permutation on each component of \Y{} so that the \msec{} is also monotone.  
		By throwing out from \Y{} any vertices or images of vertices (a finite set), we obtain
		an \msec{} satisfying \refer{eqn}{regsec2} and
		\begin{enumerate}
			\item if $y$ and \yp{} belong to the same component of \Y{}, and $1\leq j<\jp$,
			then
				$$\psisof{j}{y}\admitless\psisof{\jp}{\yp}.$$
			\item if $j\neq\jp$ and for some sequences \inY{\ys{i}, \yps{i}} we have
				$$\lim\psisof{j}{\ys{i}}=\lim\psisof{\jp}{\yps{i}}=u$$
			then $y\eqdef\fof{u}\notin\Y$.
		\end{enumerate}
	\switchtotext
		
\end{markedremark}


\subsection{Extreme Preimages}\label{subsec:extreme}

Given the linear ordering $\admitless$ on \T{} as in \refer{sec}{background},
for any point $y\in\T$ we define
	\begin{align}\label{eqn:Mymy}
		\my&\eqdef\min\fprey\\
		\My&\eqdef\max\fprey.
	\end{align}


Given the set \Y{}, we can, without reference to any \msec{}, define
	\begin{align}\label{eqn:extremes}
		\Hminus&\eqdef\clos\setbld{\my}{\inY{y}}\\
		\Hplus&\eqdef\clos\setbld{\My}{\inY{y}}\\
		\Hboth&=\Hminus\cap\Hplus.
	\end{align}

\begin{markedremark}\label{rmk:cover}
	Suppose $x\in\Hboth{}$.
	\begin{enumerate}
		\item\label{rmkcover1} If $J$ is a set which is comparable to some neighborhood 
		$U\pr$ of $x$,
		then no neighborhood $U$ of $x$ has $\fof{U}\subset\fof{J}$.
		
		\item\label{rmkcover2} In particular, $f$ does not collapse any neighborhood $U$ of $x$.
	\end{enumerate}
	\switchtotext
	To see the first statement, suppose $J\admitless U$ and $\fof{U}\subset\fof{J}$.  Then 
	for every $\xp\in U$ there exists $\xpp\in J$ with $\fof{\xp}=\fof{\xpp}\eqdef y$;  but since
	$\xpp\admitless\xp$ we have that $\xp\neq\my$;  hence $x\notin\Hs{1}$.  Similarly, if 
	$U\admitless J$ then $x\notin\Shiftset{m}$.
\end{markedremark}

\begin{lemma}\label{lem:divide}
	Suppose $x\in\Hboth$, together with \fof{x}, is not a branchpoint of \T{}. Then each 
	of the two closed intervals into which $x$ divides the edge containing it is mapped into a single
	closed branch of \T{} relative to \fof{x}.
\end{lemma}

\begin{proof}
	If not, then there exists a closed subinterval $J$ of the edge, disjoint from $x$, 
	whose image contains \fof{x} in its interior.  Then picking a neighborhood $U$ of $x$
	so small that $\fof{U}\subset\fof{J}$, we obtain a contradiction to \subref{rmk}{cover}{rmkcover1}.
\end{proof}

Given the linear ordering on \T{}, we call an \msec{} on $\Y\subset\T$ \deffont{spanning} 
if for every \inY{y}
\begin{align*}
	\psisof{1}{y}&=\my\\
	\psisof{m}{y}&=\My.
\end{align*}
\begin{prop}\label{prop:monospan}
	If \fTT{} is \comfold{}, and \admitless{} is a linear ordering on \T{} as in \refer{sec}{background},
	then there exists a monotone spanning \msec{} $\psi$ 
	defined on a cocountable open set $\Y\subset\T$ for which the corresponding $m$-shift system
	\Shiftsystem{} defined by
		$$\Shiftset{j}\eqdef\clos\setbld{\psisof{j}{y}}{\inY{y}},\quad\ilist{j}{1}{m}$$
	satisfies
		$$\fof{\kerH}\cap\Y=\emptyset,$$
	and in particular \kerH{} is eventually countable.
\end{prop}

\begin{proof}
	From \refer{prop}{regsec} and \refer{rmk}{renumber}, we can find a monotone \msec{} on a 
	cocountable open set $\Y\subset\T$ satisfying \refer{eqn}{regsec2}; using \refer{rmk}{monsec},
	we can also insure that \Y{} has no vertices or images of vertices, and that any point in \kerH{}
	has image outside \Y{}.  Now, if we replace \psisof{1}{y} \resp{\psisof{m}{y}} with \my{} \resp{\My},
	we don't change monotonicity, and we gain spanning.  We need to show that 
		\begin{align*}
			\Hminus\cap\Shiftset{j}\cap\fpreof{\Y}=\emptyset&\text{ for }j>1\\
			\Hplus\cap\Shiftset{j}\cap\fpreof{\Y}=\emptyset&\text{ for }j<m\\
			\Hboth\cap\fpreof{\Y}=\emptyset.&
		\end{align*}
	Suppose $u$ belongs to one of the sets above, and let $U$ \resp{$\tilde{U}$} 
	be a neighborhood of $u$ \resp{$y=\fof{u}$} with $\fof{U}\subset\tilde{U}$, with
	$\tilde{U}$ contained in a single component
	$\bar{\Y}$ of \Y{}, and the length of $U$ sufficiently small that any sequence \psisof{j}{\ys{i}}
	contained in $U$ has $j$ constant.  This means that there is an index
	\jp{} such that \psisof{\jp}{\yp} separates \psisof{j}{\yp} from \ms{\yp} \resp{\Ms{\yp}} for every $\yp\in U$,
	and in particular $u$ cannot be a limit of points of the form \myi{} \resp{\Myi}.
\end{proof}


\section{The Center}\label{sec:center}

In this section we concentrate on the second hypothesis of \refer{thm}{main}: the existence of a set 
containing all minimal sets in the core which locally divides the shift system \Shiftsystem{}.
Given \fTT{} a cocountably \mfold{} map on the tree \T{},   \refer{prop}{monospan} has given us a monotone spanning \msec{} $\psi$ for which the kernel is eventually countable.  Note that this implies in particular that every minimal set in the core is a periodic orbit, since infinite minimal sets are uncountable, while the core is contained in every image of the kernel.  Thus we are interested in the behavior of orbit segments near a periodic orbit in the core.

We begin by noting some simplifying assumptions concerning $f$ that we can make without loss of generality.  Note that any closed connected subset $\Tpr \subset\T$ is itself a tree.  A subtree is \deffont{natural} if it is a union of edges of \T{};  it is \deffont{proper} if it is neither all of \T{} nor a single point.

We shall concentrate on maps \fTT{} satisfying

\begin{assumption}\label{assume:assumption1}
	\begin{enumerate}
		\item\label{assumption1} If a branchpoint $y$ of \T{} is $f$-preperiodic, 
		then \fof{y} is a fixed point of $f$:
			$$\bigcup_{n=0}^{\infty}\fpre{n}{\Peri}\cap{\Branchpts}
				\subset\fpre{1}{\Fix}\cap\Branchpts.$$
		\item\label{assumption2}There is no proper $f$-invariant natural subtree of \T{}. 
	\end{enumerate}
\end{assumption}

\begin{markedremark}\label{rmk:assumptions}
	It suffices to prove \refer{thm}{main} for the class of maps satisfying \refer{assume}{assumption1}.
\switchtotext
	To see that we can assume \eqref{assumption1}, 
	note that since \Branchpts{} is a finite set, the required property holds
	for some iterate \fto{n} of $f$.  But if $f$ is cocountably \mfold{} then \fto{n} is 
	cocountably \fold{(m^{n})}, and $\ent{\fto{n}}=n\ent{f}$, so the estimates 
	$\ent{\fto{n}}\geq\log m^{n}$ and $\ent{f}\geq\log m$ are equivalent.
	
	To see that we can assume \eqref{assumption2}, 
	suppose that \Tpr{} is an $f$-invariant proper natural subtree.  We
	distinguish two subcases:
	\begin{itemize}
		\item If the restriction $f|\Tpr$ is cocountably \mfold{}, then we replace \fTT{} with 
		\selfmap{f}{\Tpr} using the fact that $\ent{f|\Tpr}\leq\ent{f}$ to complete the argument.
		\item If $f|\Tpr$ is \emph{not} cocountably \mfold{}, we collapse \Tpr{} to a point;  it is easy
		to see that the quotient space is a tree and the induced action of $f$ on this tree is again
		cocountably \mfold{}.  Since entropy is nonincreasing under factors, we are done.
	\end{itemize}
\end{markedremark}

We will assume from now on that $f$ has both properties above.

Suppose $y\in\Branchpts\cap\Fix$.  By a \deffont{branch germ} at $y$ we mean the intersection of a branch at $y$ with some (sufficiently small) neighborhood of $y$ in \T{}.  We call a branch (or edge) at $y$ \deffont{monotone} at $y$ if some branch germ maps into a single closed branch at $y$.  A \deffont{nonmonotone} branch  (or edge) at $y$ is one for which there are points arbitrarily near $y$ whose images belong to distinct branches at $y$.  

We shall call the numbering of left-open edges of \T{} (and hence the induced ordering \admitless{} on the points of \T{}) \deffont{$\boldsymbol{f}$-adjusted} if for every $y\in\Branchpts\cap\Fix$ every \emph{monotone} outgoing edge at $y$ is numbered lower than every \emph{nonmonotone} outgoing edge at $y$.  (We cannot \apriori{} rule out the possibility that the incoming branch is nonmonotone at $y$;  this must be handled separately.)

We can always pick our numbering to be $f$-adjusted, and we assume from now on that this property holds.



\subsection{Periodic Branchpoints in the Core}\label{subsec:perbranchcore}
In this subsection we show that any periodic branchpoint in the core locally divides \Shiftsystem{}.

Given $y\in\Branchpts\cap\Fix$ of outdegree \valence{}, and assuming the ordering \admitless{} is $f$-adjusted, we have numbered the branches at $y$ as
	$$\T\setminus\single{y}=\bigcup_{i=0}^{\valence}\Branchat{i}{y}$$
consistent with \admitless{}:  that is, $\Branchat{\is{1}}{y}\admitless\Branchat{\is{2}}{y}$ for $0<\is{1}<\is{2}$.   For a neighborhood \U{} of \y{} in \T{}, we use the notation
	$$\Usi\eqdef\yBranch{i}\cap\U.$$
	(Remember that \yBranch{i} does not include $y$.)
Set $\Gsys\eqdef\single{\Hone,\Hm}$ and $\Hstar\eqdef\Hone\cup\Hm$. For \ilist{i}{0}{\valence{}} and $j=1$, $m$, or $*$, let
	$$\Gij\eqdef \Usi\cap\Hs{j}.$$
	
\begin{remark}\label{rmk:collapse}
	If $\fof{\Usi}=\single{y}$ then $\Usi\cap\Hstar=\emptyset$.
\end{remark}

We will successively shrink the neighborhood \U{} to insure a number of conditions as the section progresses;  thus we are concerned with the germ of the behavior at \y{}.
In what follows, the \ith{} branch at \y{} will be denoted simply as \Branch{i}.

\begin{markedremark}\label{rmk:clint}
	Suppose $V$
	is a closed subinterval of a branch at $y$ such that 
	$\fof{V}\supsetneq\single{\y}$.  Then for a sufficiently small neighborhood \U{}
	of \y{} (in particular, one disjoint from $V$),
	\begin{enumerate}
		\item There exists at least one branch germ \Us{i} such that $\fof{V}\supset\Us{i}$.
		\item If $V\subset\Branch{\is{0}}$, where $\is{0}>0$ (\ie{} $V$ is contained in an outgoing
		branch) and  $z\in\Usi\subset\fof{V}$, then for any \ip{}, the branch germ 
		\Us{\ip} does not contain \ms{z} \resp{\Ms{z}} if $\ip> \is{0}$ 
		\resp{$\ip< \is{0}$ 
		}.
		\item If $V\subset\setbld{x}{x\admitless y}$ 
		and  $z\in\Usi\subset\fof{V}$, then for any $\ip\geq0$, the branch germ 
		\Us{\ip} does not contain \ms{z}.
	\end{enumerate}
	
	\switchtotext
	This is because if \U{} is disjoint from $V$, and then $z$ has a preimage in $V$, 
	which is \admitless{} \resp{\admitmore{}} any point of \Us{\ip} for $\ip{}\neq\is{0}$.  
	If $\ip=\is{0}$, $V$ is farther from \y{} than \Us{\ip}, which implies that 
	$V\admitless \Us{\ip}$ \resp{$V\admitmore\Us{\ip}$} when $\is{0}=0$ 
	\resp{$\is{0}>0$}.
\end{markedremark}

\begin{lemma}\label{lem:6.2}
	If \Branchi{} is not monotone at \y{}, then for \U{} sufficiently small, 
	\begin{enumerate}
		\item if $i=0$, $\Gss{i}{1}=\emptyset$;
		\item if $1\leq i\leq\valence$, $\Gss{i}{m}=\emptyset$.
	\end{enumerate}

\end{lemma}

\begin{proof}
	By \refer{rmk}{collapse}, we can assume that \Usi{} is not collapsed to $y$ by $f$.
	By assumption, there exist points $\xs{k}\in\Branchi$ converging to $y$ with $\fof{\xs{k}}=y$. 
	For each branch \Branch{\ip} intersecting \fof{\Branchi}, 
	we can find a closed interval $V_{\ip} \subset\Branchi$ such that 
	\fof{V_{\ip}} contains a neighborhood of \y{} in 
	$\single{y}\cup\Branch{\ip}$. A finite number of these suffice to fill a neighborhood
	of \y{} in \fof{\Branchi}.  Now the result follows from \refer{rmk}{clint}.
\end{proof}

Set
	$$\Phi\eqdef\setbld{i}{\Usi\cap\Hstar\nonempty
		\text{ for all neighborhoods \U{} of \y}};$$
we call branches \Branchi{} with $i\in\Ph$ \deffont{active branches} (and the corresponding
neigborhoods \Usi{} as  \deffont{active branch germs}).  
We will often refer to an active branch germ \Usi{} via just its index $i$.
We say $j\in\single{1,m}$ is a \deffont{color} for the active branch germ 
\Usi{} (or, by abuse of language, for $i$) if $\Gij\nonempty$.  Note that a given branch may have up to two colors.    A branch germ is \deffont{monochrome} if it has precisely one color.
\begin{markedremark}\label{rmk:nonmonotonecolors}
	Active nonmonotone branch germs are monochrome, with
	\begin{equation*}
		\begin{cases}
		    j=m  & \text{if $i=0$ }, \\
		    j=1  & \text{otherwise}.
		\end{cases}
	\end{equation*}

	\switchtotext
	This is an immediate corollary of \refer{lem}{6.2}.
\end{markedremark}

One can find a neighborhood \U{} of \y{} so that whenever $\tU\subset\U$ is a 
subneighborhood of \y{},  the colors of \Usi{} and \tUsi{} agree for \ilist{i}{0}{\valence{}}; we call any such \U{} a \deffont{determining} neighborhood of \y{}.
We will write 
	\begin{equation}\label{eqn:colorpass}
		\cover{i}{j}{\ip}{\jp}
	\end{equation}
if for some (hence any) determining neighborhood,
	$$\fof{\Gij}\cap\Gss{\ip}{\jp}\nonempty.$$
We form a \deffont{branch graph} \bgraph{} whose vertices are the active branch germs, and with a directed edge from \Usi{} to \Us{\ip} (denoted \bcovers{i}{\ip}) if there exist $j,\jp\in\single{1,m}$ such that \cover{i}{j}{\ip}{\jp} (in the sense of \refer{eqn}{colorpass}).

\begin{lemma}\label{lem:bgraph}
	\begin{enumerate}
		\item\label{bgraph1} Every vertex in \bgraph{} has indegree at most 2.
		\item\label{bgraph2} Every monotone vertex in \bgraph{} has outdegree at most 1.
		\item\label{bgraph3} If \bcovers{\is{1}}{\ip} and \bcovers{\is{2}}{\ip} with $\is{1}$ and $\is{2}$
			distinct, then $\is{1}$ and $\is{2}$ are both monotone and monochrome, 
			with different colors.
	\end{enumerate}
\end{lemma}

\begin{proof}
	To see \eqref{bgraph1}, suppose that \bcovers{\is{k}}{\ip} for three distinct $\is{k}$, $k=1,2,3$.
	We can assume that in an $f$-adjusted ordering 
	$\Branch{\is{1}}\admitless\Branch{\is{2}}\admitless\Branch{\is{3}}$.
	Then by \refer{rmk}{clint}, $\Gss{\is{1}}{m}=\Gss{\is{2}}{m}=\emptyset$ and
	$\Gss{\is{2}}{1}=\Gss{\is{3}}{1}=\emptyset$;  in particular, \Branch{\is{2}} is not
	active, a contradiction.
	
	\eqref{bgraph2} is clear, by the definition of monotonicity.
	
	To see \eqref{bgraph3}, suppose $\is{1}<\is{2}$, so that by \refer{rmk}{clint} 
	$\Ms{z}\not\in\Us{\is{1}}$ and $\ms{z}\not\in\Us{\is{2}}$ for any $z\in\Us{\ip}$.
	Thus we must have \Gss{\is{1}}{1} and \Gss{\is{2}}{m} both nonempty;  by 
	\refer{lem}{6.2} this means that to be active \Branch{\is{2}} must be monotone,
	and hence (since the ordering is $f$-adjusted) 
	either \Branch{\is{1}} is monotone or \Branch{\is{1}=0} is nonmonotone.
	In this last case, $\Gss{\is{1}}{1}=\emptyset$ by \refer{lem}{6.2} again, contradicting
	the assumption that \Branch{\is{1}} is active.
		 
\end{proof}

We call a path or loop in \bgraph{} \deffont{monotone} if every vertex occurring along the path is monotone at \y{}.
\begin{lemma}\label{lem:monochrome}
	Any monotone loop in \bgraph{} contains at least one monochrome vertex.
\end{lemma}

\begin{proof}
	If the loop contains all the active branch germs at \y{}, then since $f$ is cocountably \mfold{}, 
	some branch 
	contains a nontrivial preimage of \y{}, and by \refer{rmk}{clint}, this implies some 
	branch \Branchi{} has \Gss{i}{1} or \Gss{i}{m} empty.
	
	If some branch is not an element of our loop, then since the union of the (closed)
	branches in the loop is a proper natural subtree, it is not invariant (by our basic 
	assumption).  So at least one of these branches must map to a union of two or 
	more branches (one in the loop, the other out of the loop) and again 
	by \refer{rmk}{clint} the loop contains a monochrome branch germ.
\end{proof}

\begin{lemma}\label{lem:concat}
	Every monotone path \gam{} in \bgraph{} can be written as a concatenation
		$$\gamma=\concat{\alpha}{\lambda^{k}}$$
	where $\lambda$ is a loop and $\alpha$ is a path with 
	no repetitions (in particular, $\abs{\alpha}<\card\Ph$).
\end{lemma}

\begin{proof}
	This is an almost immediate consequence of the fact that every active monotone
	branch germ has outdegree 1 in \bgraph{} (\refer{lem}{bgraph}).
\end{proof}

Let \U{} be a determining neighborhood of \y{}.  We refer to a point $x$ such that $\ftoof{k}{x}\in[\U\cap\Hstar]\setminus\single{\y}$ for \ilist{k}{0}{n-1} as a \deffont{satellite} of \y{} with \deffont{time of flight} $n$.  Such a point has two kinds of itinerary: a \deffont{branch itinerary} \finseqlist{i}{k}{0}{n-1} defined by $\ftoof{k}{x}\in\Us{\is{k}}$, and a \deffont{\colorit} \finseqlist{j}{k}{0}{n-1} satisfying $\ftoof{k}{x}\in\Hs{\js{k}}$, \ilist{k}{0}{n-1}.  The branch itinerary is unique, but \apriori{} a satellite may have more than one \colorit{}.  However, \refer{lem}{6.2}, \refer{rmk}{nonmonotonecolors} and \refer{lem}{bgraph} give limitations on the \colorits{} which can occur in conjunction with a given branch itinerary.  We shall call a choice of \colorit{} \deffont{legitimate} for a given branch itinerary if it is consistent with these limitations.

\emph{Forbidden Color Words:}
We wish to find a finite color word which does not appear in any \colorit{} of any satellite of \y{}.  This 
means that \single{y} locally divides \Shiftsystem{} (\refer{dfn}{locdiv}).  We will do this in \refer{cor}{badword}, but first we need a substantial digression. 
A \deffont{colored path} of length $n$ in \bgraph{} is a path of length $n$  together with a legitimate choice of color for each branch germ \Us{\is{k}}:
	$$\gamma=(\is{0},\js{0}),\dots,(\is{n-1},\js{n-1}),
		\quad \is{k}\in\Ph,\ \js{k}\in\single{1,m}\text{ for }\ilist{k}{0}{n-1}.$$
The branch itinerary, together with a choice of  \colorit{}, for any satellite of \y{} with time of flight (at least) $n$, determines a colored path in \bgraph{} of length $n$, so the number of \colorits{} of length $n$ which occur among the satellites of \y{} is bounded above by the number 
$\boldsymbol{\Numn}$ of \colorits{} occurring among the colored paths of length $n$ in \bgraph{}. 
Since the number of (abstract) color words of length $n$ is $2^{n}$, it will suffice to prove

\begin{prop}\label{prop:Numn}
	For $n$ sufficiently large,
		$$\Numn<2^{n}.$$
\end{prop}

Given a colored path $\gamma=\ijs{0}...\ijs{n-1}$, we have a color word 
	$$\boldsymbol{\cword{\gamma}}\eqdef\js{0}...\js{n-1}$$
consisting of the sequence of colors appearing along \gam{}.  We can think of \colorword{} as a projection map from colored paths to color words.
To estimate \Numn{}, we construct, for each (legitimate) colored vertex $(b,c)\in\Ph\times\single{1,m}$, a rooted tree $\boldsymbol{\bcgraph}$ whose vertices $\boldsymbol{\gamma \in\bcn}$ at level 
\inflist{n}{0,1} are the colored paths 
	$$\gamma=\ijs{0}...[\ijs{n}=(b,c)]$$
of length $n+1$ which end at $(b,c)$, and an edge connecting each vertex $\gamma=\ijs{0}...\ijs{n}$ at level $n>0$ with the colored path $\gamp=\ijs{1}...(\is{n},\js{n})\in\bcs{n-1}$ at level $n-1$ obtained by truncating the first colored vertex.  The standard orientation of this edge is from \gamp{} to \gam{} (that is, in the direction of increasing level), which may at first appear counter-intuitive.  We shall, however, make use of this orientation only briefly (\emph{cf.} the paragraph preceding \refer{lem}{cutpt}) .  Denote the number of color words (of length $n$) occuring for vertices at level $n-1$ of \bcgraph{} by
	$$\boldsymbol{\Pnumbcn}\eqdef\card\cword{\bcs{n-1}}.$$
Clearly, $\Numn\leq\sum_{(b,c)\in\Ph\times\single{1,m}}\Pnumbc{n}$.  If $\boldsymbol{\el{}=\card\Ph}$ denotes the number of active branches at \y{}, we will establish the inequality 
\begin{equation}\label{eqn:Pnum}
	\Pnumbc{p\el+1}\leq 2^{p(\el-1)}(1+6p)
\end{equation}
for $p\geq0$ from which the proposition will follow easily.

Our previous results yield the following information about the graph \bcgraph{}:
\begin{lemma}\label{lem:bcn}
	For each vertex $\gamma=(\is{0},\js{0})\dots(b,c)$ in $\bcn$ 
	\begin{enumerate}
		\item\label{bcn1} If $n>0$ there is precisely one vertex 
			$\gamp=(\is{1},\js{1})...(b,c)\in\bcs{n-1}$ at level
			$n-1$ joined to \gam{}.
		\item\label{bcn2} In general, there are at most two vertices 
		$\gammin=(\is{-1},\js{-1})(\is{0},\js{0})...(b,c)$
			at level $n+1$ joined to \gam{}.
		\item\label{bcn3} If \Us{\is{0}} is non-monotone, then no vertex in \bcn{} other than \gam{} 
			is joined to \gamp{}. 
	\end{enumerate}
\end{lemma}

\begin{proof}
	(1) is trivial.
	
	To see (2) and (3), note that if $\gammin=(\is{-1},\js{-1})(\is{0},\js{0})...(b,c)$ 
	then \bcovers{\is{-1}}{\is{0}}.  From \subref{lem}{bgraph}{bgraph1}, there are 
	at most two possibilities for 
	\is{-1}, given \is{0};  if two distinct possibilities \iss{-1}{1}, \iss{-1}{2} exist, then both are monotone, 
	with $\jss{-1}{1}\neq\jss{-1}{2}$;  otherwise, \is{-1} is unique, and 
	can be colored in at most two ways.  By \refer{rmk}{nonmonotonecolors}, a nonmonotone germ can
	be colored in at most one way.
\end{proof}

To establish \refer{eqn}{Pnum}, we will distinguish colored paths according to the branch paths they represent.  For \inflist{p}{0}, let $\boldsymbol{\Mbcnump}$ denote the number of color words coming from colored paths of length $p\el+1$ 
in which the branch path is monotone, and $\boldsymbol{\Snumbcp}$ the number coming from colored paths of length $p\el+1$ going through at least one non-monotone branch germ.  (The reason for this peculiar numbering will become clearer in what follows.)

Let us first estimate \Mbcnump{}.  Note that $\Mbcnum{0}=1$, and if $p>0$, then for $\Mbcnump$ to be nonzero we need $b$ to belong to a monotone loop $\lambda$ in \bgraph{}.  We denote the length of $\lambda$ by $\boldsymbol{\lenlam}$; to obtain an estimate on \Mbcnump{} independent of $b$, we let $\boldsymbol{q}$ denote
the maximum length of all monotone loops in \bgraph{} (note that these are disjoint, by \subref{lem}{bgraph}{bgraph2}, and hence there are finitely many); note that 
	$$q\leq\el.$$
Given a monotone loop $\lambda$ containing $b$, \refer{lem}{monochrome} implies that 
$\lambda$ must contain at least one monochrome vertex, and hence there are at most $2^{\lenlam-1}$ legitimate colorings of $\lambda$, and at most 
	$$(2^{\lenlam-1})^{k}=2^{\len{\lambda^{k}}(1-\frac{k}{\len{\lambda^{k}}})}
		\leq 2^{\len{\lambda^{k}}\qfactor}$$
legitimate colorings of $\lambda^{k}$, the concatenation of $\lambda$ with itself $k$ times.
By \refer{lem}{concat}, every monotone path ending at $b$ is a concatenation of the form 
$\gamma=\concat{\alpha}{\lambda^{k}}$ for some $k$, where \alfa{} is a nonrepetitive monotone path, whose length is therefore bounded by the number of monotone vertices, hence by \el{}.  There is a unique path of length $p\el+1$ which is a subpath of some power of $\lambda$, and the number of colorings of it is bounded by $2^{p\el-\lfloor\frac{p\el}{q}\rfloor}\leq2^{p\el\qfactor +1}$;  there is also the possibility of replacing the initial subword of this with a nonrepetitive monotone path (\ie{} $\alpha$);  the number of legitimate colorings of \alfa{} is bounded above by $2^{\el}$.  Thus we have the estimate
	$$\Mbcnump\leq2^{p\el\qfactor+1}+2^{\el}2^{(p-1)\el\qfactor+1};$$
factoring out $2^{p\el}$ and using the fact that $q\leq\el$, we obtain the estimate
\begin{equation}\label{eqn:Mp}
	\Mbcnump\leq2^{p\el}[2^{-p}+2^{-(p-1)}]\cdot2=2^{p(\el-1)}\cdot 6.
\end{equation}

Now consider \Snumbcp{}.  To estimate this we need an excursion into abstract graph theory.  By construction, the graph \bcgraph{} is an (infinite) tree;  \refer{lem}{bcn} ((\ref{bcn1}) and (\ref{bcn2})) tells us that (if we adopt the convention that edges are oriented in the direction of increasing level) every vertex except the root  has indegree 1 and every vertex has outdegree at most 2.  We refer to such a graph as a \deffont{stump} and to a vertex with outdegree 0 \resp{1} as an \deffont{end} \resp{\deffont{cutpoint}} of the graph.  Each end has a unique path to the root; we call it a \deffont{cut end} if it is not the root, and this path contains at least one cutpoint.

\begin{lemma}\label{lem:cutpt}
	In any stump, the number of cut ends at level \el{} is bounded above by $2^{\el-1}$.
\end{lemma}

\begin{proof}
	First, assign to each edge in the tree a 0 or 1;  if the edge is leaving a cutpoint, make  sure it is 
	assigned a 0 (the assignment to edges leaving a vertex with outdegree 2 can be chosen in an
	arbitrary way).  Then each vertex at level \el{} is assigned a sequence of 0's and 1's corresponding
	to the unique path from it to the root.  If $v$ is a cut end at level $\ell$, 
	consider the sequence obtained from its
	path by replacing the first 0  associated to a cut point with a 1;  
	this leads to a sequence which does not occur in the graph,
	and is a one-to-one map from cut ends into the set of ``missing'' ends.  Since there are $2^{\el}$
	sequences all together and the set of ``missing'' sequences is disjoint from the set of extant ones,
	we have that the number of cut ends at level \el{} plus the number of ``missing'' image 
	sequences adds to at most $2^{\el}$; but the number of image sequences equals the number 
	of cut ends, and we are done.
\end{proof}

\begin{corollary}\label{cor:cuts}
	For any colored vertex $(b,c)\in\Ph\times\single{1,m}$,
	$$\Snumbc{1}\leq2^{\el-1}.$$
\end{corollary}

\begin{proof}
	Since a colored path contains at least one nonmonotone vertex, we see that some positive
	level must contain a nonmonotone vertex, and 
	it follows from \subref{lem}{bcn}{bcn3} that a nonmonotone vertex at level $k$ in \bcgraph{} 
	means that the corresponding vertex at level $k-1$ is a cut point.  In particular, the number of
	nonmonotone colored paths of length $\el+1$ ending at a given vertex is at most $2^{\el-1}$, 
	and this is a bound 
	on \Snumbc{1}. 
\end{proof}

\begin{proofof}{\refer{eqn}{Pnum}}
	If $b$ is not part of a monotone loop (in particular, if $b$ itself is not monotone), 
	then $\Mbcnump=0$ and every path ending at $(b,c)$ 
	hits a nonmonotone vertex at least once in every \el{} steps;  the number of legitimate colorings
	for all such paths is bounded above by a product of terms of the form \Snum{1}{b_{i}}{c_{i}},
	where the subscripted vertices are the ones occurring at precise multiples of \el{} steps; each of 
	these is bounded by $2^{\el-1}$, by \refer{cor}{cuts}.  
	In this case,
			$$\Pnumbc{p\el+1}\leq\Snumbcp\leq 2^{p(\el-1)}\leq 2^{p(\el-1)}(1+6p).$$
			
	When $b$ \emph{is} part of a monotone loop, we can still imagine paths ending at $(b,c)$
	of the type analyzed above, but of course there are others.  For any given path we let
	$k$ be the maximum integer for which the last $k\el+1$ vertices are monotone;  then this
	path consists of a path of length $(p-k)\el+1$ of the type above fused with a monotone path 
	of length $k\el+1$.  (The case above is $k=0$.) The first part(s) can be colored in at most 
	$2^{(p-k)(\el-1)}$ ways, as above, while the last part(s) can be colored in at most 
	\Mnum{k}{b}{c} different ways.  Using \refer{eqn}{Mp} with $k$ in place of $p$,
	we obtain \refer{eqn}{Pnum}
	\begin{align*}
		\Pnumbc{p\el+1}&\leq \sum_{k=0}^{p}2^{(p-k)(\el-1)}\Mnum{k}{b}{c}\\
			&\leq 2^{p(\el-1)}+\sum_{k=1}^{p}2^{(p-k)(\el-1)}\cdot6\cdot2^{k(\el-1)}\\
			&= 2^{p(\el-1)}[1+\sum_{k=1}^{p}6]\\
			&=2^{p(\el-1)}[1+6p]
	\end{align*}
	as required.
\end{proofof}

\begin{proofof}{\refer{prop}{Numn}}
	The number of colored vertices $(b,c)\in\Ph\times\single{1,m}$ is bounded by 
	$2\el$, so substituting
	in \refer{eqn}{Pnum} we have
	$$\Numcolors{p\el+1}\leq \sum_{(b,c)\in\Ph\times\single{1,m}}\Pnumbc{p\el+1}
		\leq 2\el\cdot2^{p(\el-1)}[1+6p].$$
	Since $1+6p$ grows more slowly than $2^{p}$, we can find a sufficiently large value
	of $p$ so that $\el\cdot[1+6p]<2^{p}$; then for $n$ equal to this value of $p\el+1$
	we have
		$$\Numn\leq2\el\cdot2^{p(\el-1)}[1+6p]<2\cdot2^{p(\el-1)}2^{p}=2^{p\el+1}=2^{n}$$
	as required.
\end{proofof}

We are now in a position to produce a forbidden color word.

\begin{corollary}\label{cor:badword}
	There exists a word in the letters \single{1,m} which does not appear in any color itinerary 
	for any satellite of \y{}.
\end{corollary}

\begin{proof}
	\refer{prop}{Numn} shows that the number of words of sufficiently long length which appear
	in legitimate colorings of paths in \bcgraph{} (which bounds the number of words appearing
	in color itineraries of satellites of \y{}) is strictly less than the number of abstract words 
	in \single{1,m}.
\end{proof}


\subsection{Periodic Non-branchpoints in the Core}\label{subsec:nonbranchper}

In the previous subsection, we showed that any periodic branchpoint in the core locally divides \Shiftsystem{}.  We now proceed to the more difficult task of finding a set containing all periodic points in the core but away from \Branchpts{} which locally divides \Shiftsystem.

\begin{definition}  We denote the set of periodic non-branchpoints in the core by
	$$\boldsymbol{\Percor}\eqdef[\centOH\cap\Peri]\setminus\Branchpts.$$  
\end{definition}

\begin{prop}\label{prop:Pclosed}
	\Percor{} is closed.
\end{prop}

\begin{proof}
	We begin with a few observations:
	
	\begin{claim}
		A periodic point which is an accumulation point of \Percor{} 
		has an orbit disjoint from \Branchpts{}.
	\end{claim}
	This is an immediate corollary of \refer{cor}{badword}: since \centOH{} is closed, the orbit 
	belongs to \centOH{};  but then if it intersects \Branchpts{}, by our assumptions on $f$ it
	consists of a fixedpoint in \Branchpts{}, and since it is an accumulation point of \Percor{},
	there exist periodic points in \centOH{} arbitrarily near to (but distinct from) the point;  
	they are satellites of the fixedpoint with arbitrarily long time of flight, and  since they
	belong to \centH{}, they have all possible words in their itinerary, contrary to \refer{cor}{badword}.
	\subqed{}

	Now, suppose \single{\qs{n}} is a sequence of points of \Percor{} converging to $y\not\in\Percor$;
	denote by \Qs{n} the orbit of \qs{n}.
	
	Since $y\in\clos\Percor\subset\centOH$, its $\omega$-limit set must also belong to the invariant 
	closed set \centOH{}, and since the latter is contained in the countable set \ftoof{i}{\kerH}, 
	any minimal subset is  a periodic orbit.  Thus we can pick a cycle
	$P=\single{\psub{0},...,\psub{N-1}}$ 
	in $\omega(y)$, 
	which \emph{a fortiori} also belongs to \Percor{}.  We can pick a neighborhood $U$ of $P$ 
	which is disjoint 
	from  \Branchpts{} (by the claim) and from some neighborhood of $y$; 
	going to a subsequence if necessary, 
	we can also pick points $\qps{n}\in\Qs{n}$ converging to \psub{0} from one side.  
	
	Now pick a closed one-sided neighborhood $\Js{0}\subset U$ of \ps{0} containing \qps{n} for 
	all sufficiently large $n$.  We know from \refer{lem}{divide} and the fact that the \qps{n} are 
	periodic points whose (forward) orbit leaves $U$ that $\Js{1}\eqdef\fof{\Js{0}}$ 
	is again a closed one-sided neighborhood of \psub{1};  iterating this procedure 
	(reducing \Js{0} if necessary) we obtain closed one-sided neighborhoods
	$\Js{i}\eqdef\ftoof{i}{\Js{0}}\subset U$ of \psub{i}, \ilist{i}{0}{N-1}. 
	Note that \fof{\Js{N-1}} is a closed one-sided neighborhood of \psub{0}, 
	on either the same or the opposite side 
	as \Js{0}; in the first case, it must properly contain \Js{0}  (otherwise the union 
	$\bigcup_{i=0}^{N-1}\Js{i}$ is an invariant neighborhood of $P$ in $U$, contradicting the fact that
	$\qs{n}\to y$) while in the second case we can iterate the procedure up to \Js{2N}, which
	properly contains either \Js{0} or \Js{N};  in any case, (again reducing \Js{0} if necessary) we
	obtain a family of $N$ disjoint closed intervals $\tJs{i}=\Js{i}$ or $\Js{i}\cup\ftoof{N}{\Js{i}}$, 
	\ilist{i}{0}{N-1} contained in $U$ 
	with $\psub{i}\in\tJs{i}$, 
	(possibly as an endpoint), such that $\fof{\tJs{i}}=\tJs{i+1}$ for \ilist{i}{0}{N-2}, and 
	$\tJs{0}\subsetneq\fof{\tJs{N-1}}\subset U$.
	
	Since the periodic orbits \Qs{n} intersect both a neighborhood of $y$ and 
	$V\eqdef\bigcup_{i=0}^{N-1}\tJs{i}$, for each sufficiently large $n$ we can find 
	$q=\qpps{n}\in\Qs{n}$ such
	that $q\not\in V$ but $\fof{q}\in \interior V$, say $\fof{q}\in\interior\tJs{i+1}$.  
	We know that some neighborhood $U\pr$ of $q$ is comparable to \tJs{i}.  
	But since \fof{\tJs{i}} contains \tJs{i+1}, it also
	contains the image of a neighborhood of $q$, and this contradicts the assumption that 
	$q=\qpps{n}\in\Hboth$,
	by \subref{rmk}{cover}{rmkcover1}, since $\tJs{i}\subset U$ contains no branchpoints 
	and $q\in\Percor{}$.
\end{proof}

\begin{lemma}\label{lem:order}
	Suppose \psub{i}, $i=1,2,3$ are distinct elements of \Percor{} contained in a single edge, with \psub{2}
	between \psub{1} and \psub{3}.  If \fof{\psub{1}} and \fof{\psub{3}} lie in the same edge, then \fof{\psub{2}} lies
	between them;  in particular, it belongs to the same edge.
\end{lemma}
\begin{proof}
	Since they are distinct periodic points, their images under $f$ are also distinct;  \refer{lem}{divide} 
	applied to \psub{1} \resp{\psub{3}} implies \fof{\psub{1}} \resp{\fof{\psub{3}}} cannot separate \fof{\psub{2}} 
	from \fof{\psub{3}} \resp{\fof{\psub{1}}}, and if \fof{\psub{1}} and \fof{\psub{3}} lie in the same edge, this implies
	\fof{\psub{2}} lies between them on this edge.
\end{proof}

\begin{definition}\label{dfn:edgeitin}
	The \deffont{edge itinerary} of a point $p\in\Percor$ of least period $N$ is the sequence 
	$\edgeitin{p}=\single{\Es{0},...,\Es{N-1}}$ of edges of \T{} visited by $p$ during one period:
	$\ftoof{i}{p}\in\Es{i}$, \ilist{i}{0}{N-1}.  
	
	A \deffont{$\boldsymbol{t}$-fold repetition} of \single{\Es{0},...,\Es{N-1}} 
	(\inNat{t}) is a sequence of edges
	\single{\Eps{0},...,\Eps{tN-1}} with $\Eps{i}=\Es{j}$ whenever $i\equiv j\mod{N}$;  in particular, it is a 
	\deffont{doubling} if $t=2$.  An edge itinerary is \deffont{repetitive} if it is a \fold{t} repetition 
	of some shorter sequence of edges.
	
	We call two points $p,q\in\Percor$ \deffont{edge equivalent} (and write \edgequiv{p}{q}) if one
	of their edge itineraries is a \fold{t} repetition of the other for some \inNat{t}.  Denote the 
	edge equivalence class of $p\in\Percor$ by $\boldsymbol{\edgeclass{p}}$, and its convex 
	hull in \T{} by $\boldsymbol{\edgehull{p}}$.
\end{definition}

\begin{prop}\label{prop:equiv}
	\begin{enumerate}
		\item\label{equiv1} If $p\in\Percor$ and $t>2$ then \edgeitin{p} cannot be a \fold{t} 
		repetition of any sequence of edges.  Thus at most two itineraries can occur among
		the elements of one edge equivalence class, and in this case one is a doubling of 
		the other.
		
		\item\label{equiv2} If \edgequiv{p}{q} with \edgeitin{q} doubling \edgeitin{p},
		then $p$ lies between $q$ and 
		\edgequiv{\ftoof{N}{q}}{q},
		where $N$ is the period of $p$.
		
		\item\label{equiv3} For any $p\in\Percor$, 
			\begin{enumerate}
				\item \edgehull{p} is a closed (possibly degenerate) interval 
				interior to a single edge;
				
				\item  \edgehull{p} and \edgehull{q} are disjoint unless \edgequiv{p}{q};
				
				\item $\fof{\edgehull{p}}=\edgehull{\fof{p}}$.
			\end{enumerate}
	\end{enumerate}

\end{prop}

\begin{proof}
	\begin{enumerate}
		\item If for some $p\in\Percor$ the itinerary \edgeitin{p} is a $t$-fold repetition of an 
		itinerary of length $k$, then its period $N=kt$ and the orbit of $p$ under \fto{k} 
		consists of $t$ points,
		all belonging to \Percor{} and having the same edge itinerary, which are permuted by \fto{k}.
		Inductive application of \refer{lem}{order} shows that the extreme points of this set are
		preserved by \fto{k}, which means that if $t\geq3$, no internal point of this set can map to
		an extreme point of its image under any iterate of $f$, 
		contradicting the assumption that they form a periodic orbit under $f$.
		
		\item Suppose \edgequiv{p}{q} and the period of $q$ is greater than that of $p$; 
		denote the period of $p$ by $N$.
		By the preceding
		item, \edgeitin{q} is a doubling of \edgeitin{p}, so $S\eqdef\single{p,q,\ftoof{N}{q}}$ 
		is contained in 
		one edge, for each $i$, \ftoof{i}{S} is contained in a single edge,
		and $\ftoof{N}{S}=S$.
		By \refer{lem}{order}, the middle point of $S$ cannot map to either extreme point
		of $S$, which forces $p$ to be the middle point of $S$.
		
		\item \begin{enumerate}
				\item This is a trivial consequence of (\ref{equiv1}) and \refer{prop}{Pclosed}.
				
				\item We need to show that if $p\in\Percor$ belongs to \hull{q}
				then \edgequiv{p}{q}.  Suppose $p\in\hull{q}$ but \nedgequiv{p}{q}.
				Then there exists \edgequiv{\qpp}{q} so that $p$ lies between
				$q$ and \qpp{}, and by \refer{lem}{order} the same holds for 
				the images of these three points by any iterate \fto{i}.
				But then \refer{lem}{order} forces $q$ to have the same 
				edge itinerary as $p$, or else
				to be a doubling of it.
				
				\item We need to show that a point $x\in\edgehull{p}$ 
				(not assuming $x\in\Percor$) has $\fof{x}\in\edgehull{\fof{p}}$.  
				This follows from application of 
				\refer	{lem}{divide} to the endpoints of \edgehull{p}, 
				together with the observation that
				\edgehull{p} contains no \bpt{}s, by the earlier arguments in this proof.

			\end{enumerate}
		
			\end{enumerate}

\end{proof}

\begin{lemma}\label{lem:finite}
	\Percor{} has finitely many edge-equivalence classes.
\end{lemma}

\begin{proof}
	Suppose $\psub{n}\in\Percor$, \inflist{n}{1,2} belong to pairwise non-edge-equivalent orbits, 
	with respective periods
	\Ns{n}.  Passing to a subsequence, we can assume the \psub{n} converge inside a single edge 
	$\Es{0}$ to a point $p$ which by \refer{prop}{Pclosed} belongs to \Percor{} and hence is periodic,
	say with period $N$.  For each $n$, let 
		$$\ts{n}\eqdef\min
		\setbld{t>0}{\ftoof{\Ns{n}-t}{\psub{n}},\ftoof{N-t}{p}\text{ belong to different edges}}.$$
	We can assume by passing to a subsequence that all the \ts{n} are congruent modulo $N$,
	and all the \ftoof{\Ns{n}-\ts{n}}{\psub{n}} belong to the same edge.
	By applying an iterate of $f$ to the whole picture, we can assume that $\ts{n}=1$ for all $n$:  thus
	we have a sequence $\qs{n}\eqdef\ftoof{\Ns{n}-1}{\psub{n}}$ all contained in an edge
	different from that containing $\pp\eqdef\ftoof{N-1}{p}$.
	
	Let $q$ be an accumulation point of the \qs{n}; by \refer{prop}{Pclosed} $q\in\Percor$ is 
	periodic, and distinct from \pp{}; but both are periodic, and both map to $p$, a contradiction.
\end{proof}

For any $p\in\Percor$, the union
	$$\Zp\eqdef\bigcup\ftoof{i}{\edgehull{p}}$$
consists of $N$ disjoint closed (possibly degenerate) intervals
	$$\Zp=\bigcup_{j=0}^{N-1}\Zs{j}$$
where $N$ is the least period among the points in the edge equivalence class \edgeclass{p} and the numbering is via the action of $f$:
	$$\fof{\Zs{j}}=\Zs{j+1}$$
(indices taken mod $N$).

\begin{markedremark}\label{rmk:Zmfold}
	If for some $p\in\Percor$ the restriction $f|\edgehullorb{p}$ is cocountably \mfold{}, then 
		$$\ent{f}\geq\log m.$$
	\switchtotext
	This is Theorem 4.3 in \cite{Bobok:mfold}, which gives our desired inequality for \mfold{}
	interval maps,  applied to $\fto{N}|\Zs{0}$, which is 
	nondegenerate and cocountably \fold{(m^{N})}.
\end{markedremark}
In view of \refer{rmk}{Zmfold}, we can \emph{assume for the rest of this section} that the following holds:
\begin{assumption}\label{assume:noncomfold}
	The restriction of $f$ to each set \edgehullorb{p}, $p\in\Percor$, fails to be cocountably \mfold{}.
\end{assumption}

If \Zp{} is not connected ($N\geq2$), at least one component of the complement \Zcompl{} has common boundary points with at least two of the intervals \Zs{j}.  We call the closure of such a component a 
\deffont{central component} relative to \Zp{}, and denote the union of all the central components by 
$\boldsymbol{\Centcomp}$.  A \deffont{peripheral component} of \Zcompl{} is one attached to a unique interval \Zs{j}; we label its closure \Psp{j} and denote the union of these by $\boldsymbol{\Percomp}$.  This yields a partition of the tree corresponding to any 
edge equivalence class \edgeclass{p}
	$$\T=\Centcomp\cup\Zp\cup\Percomp.$$
Each interval \Zs{j} of \Zp{} touches \Centcomp{} in at least one endpoint;  we will write
	$$\Zs{j}=\clint{\zmins{j}}{\zpls{j}}$$
where \zmins{j} is a common boundary point of \Zs{j} and \Centcomp{};  if \Zs{j} is nondegenerate, then the other endpoint \zpls{j} may touch another central component, or the peripheral component \Psp{j}, or be an end of \T{}, in which case we shall refer to a ``trivial peripheral component'' $\Psp{j}=\emptyset$.
When \Zs{j} is a single point, we shall nonetheless refer separately to its endpoints \zmins{j} and \zpls{j}.

If $N=1$, \ie{} $\Zp=\Zs{0}=\clint{\zmins{0}}{\zpls{0}}$, then \Zcompl{} has at most two components, both peripheral;  we shall associate to each endpoint a (possibly trivial) peripheral component \Pspeither{0}.  In this case, our \refer{assume}{noncomfold} says that at least one of the peripheral components has interior points mapping to \Zp{}, and we make sure \Pspmax{} is of this type.

Under \refer{assume}{noncomfold}, we will construct a set \Sepset{} containing \Percor{} which locally divides \Shiftsystem{}.  This will be made up of sets of the form $\sepset{p}\supset\Zp$ for various edge equivalence classes in \Percor.  Each such set \sepset{p} will be either \Zp{} itself or its union with \Centcomp{} or \Percomp{} if either contains no interior points mapping to \Zp{}.  We begin by showing that if this occurs then the partition $\T=\Centcomp\cup\Zp\cup\Percomp$ has a particular structure.

Suppose first that $\fpre{1}{\Zp}\cap\Percomp=\emptyset$.

\begin{markedremark}\label{rmk:perinvar}
	If $N\geq2$ and no point interior to \Percomp{} maps into \Zp{}, then either $\Percomp=\emptyset$
	(so \zpls{i} is an endpoint of \T{} for \ilist{i}{0}{N-1}) or else
	each \Zs{i} is attached to a nontrivial peripheral component \Psp{i}, and for \ilist{i}{0}{N-1}
	$\fof{\Psp{i}}\subset\Psp{i+1}$ (indices taken mod $N$).  In either case, \Zcompl{} has exactly one 
	central component.
\switchtotext
	This follows immediately from \refer{lem}{divide}.
\end{markedremark}

We would like to establish an analogous picture when no interior point of \Centcomp{} maps to \Zp.

\begin{prop}\label{prop:onecent}
	If $\Centcomp\neq\emptyset$ but $\fpre{1}{\Zp}\cap\Centcomp=\emptyset$, then
	\begin{enumerate} 
		\item \Zcompl{} has a unique central component, 
		\item each component \Zs{i} is attached to a nontrivial peripheral component \Psp{i};
		\item for \ilist{i}{0}{N-1} (taking indices mod $N$),
		$$\Psp{i+1}\subset\fof{\Psp{i}}\subset\Centcomp\cup\Zs{i+1}\cup\Psp{i+1}.$$
	\end{enumerate}
\end{prop}

\begin{proof}
	To establish the first statement, we will show that if \Zcompl{} has 
	at least two central components, then some interior point of \Centcomp{} 
	maps into \Zp{}.
	If there are at least two central components, \Zp{} includes three intervals 
	$\Zis{j}$, \ilist{j}{1,2}{3}, such that $\Zis{2}\subset\hull{\Zis{1},\Zis{3}}$,
	and no other intervals \Zs{i}, $i\neq\is{1},\is{2},\is{3}$ intersect \hull{\Zis{1},\Zis{3}}.
	 Since $\fof{\Zs{i}}=\Zs{i+1}$ for all $i$ (indices taken mod $N$) and the \Zs{i} are 
	 permuted, we can assume that \Zs{\is{2}-1} is not contained in \hull{\Zs{\is{1}-1},\Zs{\is{3}-1}}.
	 But then the interior of this hull is disjoint from all peripheral components and contains a
	 point mapping into \Zis{2}, which must then be interior to \Centcomp{}.

	This establishes the uniqueness of the central component.  
	If \Zs{i} touches a peripheral component, we have already called it \Psp{i};
	if not, then \zpls{i} is an endpoint of the tree \T{}, and we say \Psp{i} is trivial.
	Since $f|\Zp$ is not cocountably \mfold{} and \Zp{} has no preimages in \Centcomp{},
	at least one peripheral component contains
	an interior point mapping into \Zp{}, and in particular at least one \Psp{i} is nontrivial.  
	
	Now, suppose that for some $j$, \Psp{j} is nontrivial.
	Since the map $f$ is surjective on \T{} and $\Centcomp\cup\Zp$ is invariant, 
	$\Psp{j}\subset\fof{\Percomp}$.  We claim that only \Psp{j-1} (mod $N$) can contain preimages
	of interior points of \Psp{j}.  Note that since $\Centcomp\neq\emptyset$ requires $N\geq2$, the
	points \zpls{i} are not fixedpoints, and in particular (by \refer{assume}{assumption1}) 
	have no preimages which are branchpoints.
	Now suppose some $\zp\in\Psp{i}$, $i\neq j-1$, has $\fof{\zp}\in\interior\Psp{j}$.  Then
	the interval \clint{\zpls{i}}{\zp} maps across a neighborhood of \zpls{j}, which we can assume
	to contain no branchpoints.  By \refer{rmk}{intervals}, we can find a subinterval $J$ of 
	\clint{\zpls{i}}{\zp} which is disjoint from \Branchpts{} mapping exactly onto this neighborhood;
	but then some neighborhood $U$ of \zpls{j-1} (also disjoint from \Branchpts{}) has 
	$\fof{U}\subset\fof{J}$;  since both $U$ and $J$ are comparable and disjoint from \Branchpts{},
	\subref{rmk}{cover}{rmkcover1} tells us that $\zpls{j}\notin\Hboth$, 
	contradicting $\zpls{j}\in\centOH$.	
	
	Note, however, that this shows that if \Psp{j} is nontrivial, then so is \Psp{j-1}, and no other
	peripheral branch can have points mapping to its interior.  Inductively, this proves the proposition.
\end{proof}

To continue our analysis, we need to track the dynamics of neighborhoods of the points \zpms{j}.
To this end, we assign to each endpoint \zpms{j} of \Zs{j} a set \Upms{j} as follows:
\begin{enumerate}
	\item If \zpls{j} is an endpoint of \T{}, then $\Ups{j}\eqdef\Zs{j}$;
	\item Otherwise, \Upms{j} is a one-sided neighborhood of \zpms{j}, contained in the 
	closed component of \Zcompl{} attached to \zpms{j}.
\end{enumerate}
Reducing the sets \Upms{j} in case (2) if necessary, we can assume that they are pairwise disjoint,
and that the interior of each one is disjoint from the finite set $\Branchpts\cup\fof{\Branchpts}$.

Given $j\in\single{0,...,N-1}$ and $\sigma\in\single{+,-}$, we can track the action of $f$ on \Uss{\sigma}{j}
by specifying a set $\vphiof{\Uss{\sigma}{j}}\in\single{\Ums{j+1},\Ups{j+1},\Zs{j+1}}$, well defined in
view of \refer{lem}{divide} applied to $\zpms{j}\in\Percor{}\subset\T\setminus\Branchpts$, according to
\begin{enumerate}
	\item If $\fof{\Uss{\sigma}{j}}\subset\Zs{j+1}$, then $\vphiof{\Uss{\sigma}{j}}\eqdef \Zs{j+1}$
	\item Otherwise, $\vphiof{\Uss{\sigma}{j}}\neq\Zs{j+1}$ has nontrivial intersection with
	$\fof{\Uss{\sigma}{j}}$.
\end{enumerate}
Extending $\varphi$ to $\vphiof{\Zs{j}}=\Zs{j+1}$, we have a self-map of the finite set 
\setbld{\Ums{j},\Ups{j},\Zs{j}}{j=0,...,N-1} into itself, each orbit of which is eventually periodic.
 Automatically, $\varphi$ has the \deffont{trivial} cycle 
 $\Zs{0}\mapsto\Zs{1}\mapsto...\mapsto\Zs{N-1}\mapsto\Zs{0}$.  As to nontrivial cycles, the two possibilities 
 are:
 \begin{enumerate}
	\item A single nontrivial cycle of length $2N$, with $\vphitoof{N}{\Upms{j}}=\Uss{\mp}{j}$ for each $j$
	(and all $\Uss{\pm}{j}\neq\Zs{j}$);
	\item one or two nontrivial cycles (disjoint if there are two) of length $N$ 
	($\vphitoof{N}{\Uss{\pm}{j}}=\Uss{\pm}{j}$).
\end{enumerate}

Let $U\eqdef\bigcup_{j=0}^{N-1}\Ums{j}\cup\Ups{j}$;  abusing the terminology of 
\S\ref{subsec:perbranchcore}, we refer to a point $x\in U$ as a \deffont{satellite} of \Zp{} with \deffont{time of flight} $t$ if $\ftoof{k}{x}\in U\setminus\Zp$ for \ilist{k}{0}{t-1}.  If $t\geq2N$, it is a 
\deffont{persistent satellite}.  The \deffont{local itinerary} of a satellite $x$ is the sequence of nontrivial sets \Uss{\sigma}{j} containing successive iterates of $x$:  note that if $\ftoof{k}{x}\in\Uss{\sigma}{j}\setminus\Zp$ and $k<t-1$ 
then $\ftoof{k+1}{x}\in\vphiof{\Uss{\sigma}{j}}\setminus\Zp$.  We are interested in satellites whose orbit segment is contained in $\Hstar\eqdef\Hone\cup\Hm$;  for any such satellite, a color itinerary is a word $\ws{0},...,\ws{t-1}$ with $\ws{i}\in\single{1,m}$ such that
$\ftoof{i}{x}\in\Hs{\ws{i}}$ for \ilist{i}{0}{t-1}.

To build a set \sepset{p} that locally divides \Shiftsystem{} we will try to find words in the letters \single{1,m} which 
do not occur in any color itinerary of any persistent satellite of \Zp{}.

\begin{remark}\label{rmk:trivorb}
	\begin{enumerate}
		\item If $\vphitoof{N}{\Uss{\sigma}{j}}=\Zs{j}$, then \Uss{\sigma}{j} contains no 
		persistent satellites of \Zp{}.
		
		\item Any persistent satellite $x$ has \ftoof{N}{x} contained in an element of a
		nontrivial $\varphi$-cycle.
	\end{enumerate}

\end{remark}

The following observations lead to a lemma which will be useful in proving the existence of ``forbidden words'' for satellites of 
\Zp{}.  We will say that a point \zpp{} \deffont{separates} the points $z$ and \zp{} if
it belongs to the interior of their convex hull (this means that $z$ and \zp{} belong to different branches at \zpp{}).

\begin{lemma}\label{lem:doublecover}
	Suppose $y\in\hull{\zp,z}$ such that $\interior\fof{\hull{y,z}}$ is disjoint from 
	\hull{\fof{\zp},\fof{z}}, \Branchpts{}, and \fof{\Branchpts}. 
%
	Then $\interior\hull{y,z}$ is disjoint from at least one of the sets \Hone{}, \Hm{}.
\end{lemma}

\begin{proof}
	First apply \refer{rmk}{intervals} to the interval $\interior\fof{\hull{y,z}}$ to find a subinterval
	of \hull{y,z} mapping exactly onto it;  we can assume without loss of generality that $y$ is an 
	endpoint mapping to an endpoint of \fof{\hull{y,z}}, 
	so that $\fof{\hull{y,z}}=\hull{\fof{y},\fof{z}}$.  Now apply \refer{rmk}{intervals} again,
	this time to find a subinterval $J\pr$ of \hull{\zp,y} mapping onto $\interior\fof{\hull{y,z}}$.
	By assumption, $J\pr$ is disjoint from \Branchpts{} and hence comparable to \hull{y,z}, and the 
	lemma follows from \subref{rmk}{cover}{rmkcover1}.
\end{proof}

The next three lemmas provide the basis for constructing ``forbidden words'' for satellites of \Zp{}.
We fix $p\in\Percor{}$ for these three results.

\begin{lemma}\label{lem:switch}
	Suppose $\interior\Ums{j}$ is contained in a central component of \Zcompl{} and 
	$\vphiof{\Ums{j}}\subset\Ps{j+1}$.  Then at least one of $\Hone\cap\interior\Ums{j}$
	and $\Hm\cap\interior\Ums{j}$ is empty.
\end{lemma}

\begin{proof}
	Let $z=\zmins{j}$, $y$ an endpoint of \Ums{j}, and pick 
	$\zp\in\Zp\setminus\Zs{j}$ another point in the 
	central component of \Zcompl{} containing \Ums{j}.  Since 
	$\fof{\zp}\in\Zp\setminus\Zs{j+1}$ we have the situation of \refer{lem}{doublecover}
	and our conclusion follows.
\end{proof}

\begin{lemma}\label{lem:peripheral}
	Suppose \single{\Ups{0},...,\Ups{N-1}} is a nontrivial $\varphi$-cycle with $\Ups{i}\subset\Ps{i}$
	for \ilist{i}{0}{N-1}.  If there exist points of $\Percomp\setminus\Zp$
	mapping into \Zp{}, then for some $i\in\single{0,...,N-1}$ and $j\in\single{1,m}$,
	$\Hs{j}\cap\Ups{i}=\emptyset$.
\end{lemma}

\begin{proof}
	Suppose $\zp\in\Ps{i}\setminus\Zp{}$ with $\fof{\zp}\in\Zs{i+1}$, and let
	$\Ups{i}=\ropint{\zpls{i}}{\zpls{i}+\eps}$.  Then \refer{lem}{doublecover} applied to 
	$z=\zpls{i}$, $y=\zpls{i}+\eps$ and $\zp$ gives the desired conclusion.
\end{proof}

\begin{lemma}\label{lem:central}
	Suppose a nontrivial $\varphi$-cycle has all elements in central components of \Zcompl{}.
	If some point $\zp\not\in\Zp$ contained in a central component maps into \Zp{}, then some
	element of the cycle has interior disjoint from either \Hone{} or \Hm{}.
\end{lemma}

\begin{proof}
	Suppose $\fof{\zp}\in\Zs{j+1}$, and $\Ums{j}=\hull{\zmins{j}-\eps,\zmins{j}}$; apply 
	\refer{lem}{doublecover} to $z=\zmins{j}$, $y=\zmins{j}-\eps$ and \zp{}.
\end{proof}

\begin{prop}\label{prop:Wp}
	Suppose $p\in\Percor$ with $f|\Zp$ not cocountably \mfold{}.  Define a set $\sepset{p}\supset\Zp$
	as follows:
	\begin{enumerate}
		\item If \Zcompl{} has exactly one central component, and no interior point of
			\Centcomp{}  \resp{of \Percomp} maps to \Zp{}, set
			$$\sepset{p}\eqdef\Centcomp\cup\Zp\quad\text{\resp{$\Percomp\cup\Zp$}}$$
		Note that in these cases $N\geq2$.
		\item If $N=1$ and one peripheral component \Pspmin{} has no interior points
			mapping to \Zp{} or to the other peripheral component, then
				$$\sepset{p}\eqdef\Zp\cup\Pspmin{}.$$
		\item In all other cases, 
			$$\sepset{p}\eqdef\Zp.$$
	\end{enumerate}
	
	Then 
	\begin{enumerate}
		\item $\fof{\sepset{p}}\subset\sepset{p}$.
		\item For any open set $U$ containing \Zp{}, $U\setminus\sepset{p}\neq\emptyset$.
		\item \sepset{p} locally divides \Shiftsystem{}.
	\end{enumerate}
\end{prop}

\begin{proof}
	The first property is trivial, and the second is nearly so: in any case, since some points outside
	\Zp{} map to it, \Zcompl{} is nonempty;  furthermore, when we do adjoin a set to \Zp{} 
	to form \sepset{p}, there is always another endpoint of each \Zs{j} to which a nontrivial component
	of \Zcompl{} is attached.
	
	We are left with the third property.  We need to show that every persistent satellite of \Zp{} has 
	its \Gsys-itinerary in some proper subset $\Lambda(p)\subset\mshiftspace$ 
	which depends only on \Zp{}.
	
	Since every persistent satellite lands in a nontrivial $\varphi$-cycle after at most $N$ applications
	of $f$, we limit our attention to these.  
	
	If a nontrivial $\varphi$-cycle includes subsets of both 
	central and peripheral components of \Zcompl{}, then \refer{lem}{switch} insures that a persistent 
	satellite with local itinerary in this cycle cannot have a color itinerary containing $N$ successive
	occurences of some $j\in\single{1,m}$ (where $j$ depends only on the cycle).  
	Note that 
	this situation \emph{must} occur if either the cycle has length $2N$ or the cycle contains 
	at least one peripheral element and the total number of central components exceeds 1.
	
	If the cycle is entirely peripheral \resp{entirely central} and the peripheral \resp{central} 
	components of \Zcompl{} contain at least one preimage of \Zp{},
	then \refer{lem}{peripheral} \resp{\refer{lem}{central}} again insures that any satellite with
	local itinerary in this cycle cannot display $N$ successive occurences of some $j\in\single{1,m}$.
	
	The exception that remains is when the union of the peripheral \resp{central} components of 
	\Zcompl{} is $f$-invariant and contains our $\varphi$-cycle.  In these cases, we know that the
	union \sepset{p}{} of \Zp{} with all the peripheral \resp{central} components is
	$f$-invariant, and that there must be preimages of \Zp{} outside \sepset{p}{}.  
	In particular, no open set containing \Zp{} is contained in \sepset{p}.
	If there is no nontrivial $\varphi$-cycle outside \sepset{p}{}, then \sepset{p}{} has  
	no persistent
	satellites, while if there is one then \refer{lem}{peripheral} \resp{\refer{lem}{central}} applies to it.
\end{proof}

\begin{prop}\label{prop:locsepset}
	Under \refer{assume}{noncomfold}, we can find a subset $\calQ\subset\Percor$ such that
		$$\sepset{\calQ}\eqdef\bigcup_{q\in\calQ}\sepset{q}$$
	is a proper subset of \T{} which contains $\edgehullorb{\Percor}\eqdef\bigcup_{p\in\Percor}\Zp$.
\end{prop}

\begin{proof}
	Note that if, for some $p\in\Percor$, \sepset{p} is a connected proper superset of \Zp{}, then either
	\begin{enumerate}
		\item\label{cond1} $\sepset{p}=\Centcomp\cup\Zp$, or
		\item\label{cond2} \Zp{} has a single component \Zs{0}, and $\sepset{p}=\Zs{0}\cup\Pspmin$.
	\end{enumerate}
	If any such points $p\in\Percor$ exist, we pick \psub{0} if possible of the first type and in any case so
	that \sepset{\psub{0}} is maximal (in the sense that it is not a proper subset of \sepset{p} for any
	$p\in\Percor$ of the same type).  Now set
		$$\calQ\eqdef\single{\psub{0}}\cup\setbld{q\in\Percor}{\Zq\not\subset\sepset{\psub{0}}}.$$
	In this case, we 
	\begin{claim}
		For any $q\in\calQ\setminus\single{\psub{0}}$, $\sepset{q}\cap\sepset{\psub{0}}=\emptyset$.
	\end{claim}
	
	To see this in case \ref{cond1}, note that every such \Zq{} must be contained in 
	peripheral components of \Zcomplpo{}, of which there is more than one and by 
	\refer{prop}{onecent}, these are permuted transitively by $f$.  In particular, $q$ cannot
	also be of type \ref{cond1}, since then its central component $C(q)$ (and hence \sepset{q}) 
	would contain \sepset{\psub{0}}, contradicting the maximality of the latter.  But then \sepset{q} can 
	consist only of \Zq{} possibly together with peripheral components of \Zcomplq{}, which are 
	separated from \sepset{\psub{0}} by \Zq{}.
	
	In case \eqref{cond1} fails but \eqref{cond2} holds, each \Zq{} is contained in the single peripheral
	component \Pspmaxpo{}, and \sepset{q} consists of \Zq{} with the possible addition of 
	peripheral components of \Zcomplq{}, which again are separated from \sepset{\psub{0}}
	by \Zq{}.\subqed{}
	
	Now, if no $p\in\Percor$ satisfies \eqref{cond1} or \eqref{cond2}, then for each $p\in\Percor$
	either \sepset{p} consists of \Zp{} together with peripheral components, or else
	$\sepset{p}=\Zp{}$.   If \sepset{p}=\Zp{} for all $p\in\Percor$, then clearly $\calQ=\Percor$
	has the required property, while if some $p\in\Percor$ has $\sepset{p}=\Zp\cup\Percomp$
	(with $\Percomp\neq\emptyset$), then \Centcomp{} is connected and we can pick \psub{0}
	of this type for which \sepset{\psub{0}} is maximal.  Then 
		$$\calQ=\single{\psub{0}}\cup\setbld{q\in\Percor}{\Zq\subset C(\psub{0})}$$
	works, since any $q\in\calQ\setminus\single{\psub{0}}$ with $\sepset{q}\neq\Zq$ has
	$P(q)$ separated from \sepset{\psub{0}} by \Zq{}.
\end{proof}

\begin{markedremark}\label{rmk:locsepset}
	If \calQ{} is a set of the type described in \refer{prop}{locsepset}, then
		$$\Sepset\eqdef\sepset{\calQ}$$
	contains \Percor{} and locally divides \Shiftsystem{}.
	\switchtotext
	To see this, note that since $\fof{\Shiftset{i}}=\T$ for all $i$ and $\Sepset\neq\T$, it follows
	that $\Shiftset{i}\setminus\Sepset\neq\emptyset$ for each $i$, and \refer{prop}{Wp} and
	\refer{lem}{unionlocdiv} insures
	that every orbit segment in a neighborhood of \Sepset{} either terminates in \Sepset{} or 
	has its itineraries in a set $\Lambda\subset\mshiftspace$ with entropy strictly less than $\log m$.
\end{markedremark}

\section{Proof of Main Theorem}\label{sec:main}
 In this section we indicate how the various threads of this paper can be pulled together to prove \refer{thm}{main}, by fulfilling the hypotheses of \refer{thm}{mshiftent}.  The key point is to combine  \refer{prop}{monospan}, giving a spanning \msec{} with eventually countable kernel, and  the results of \refer{sec}{center}, characterizing the center of a spanning section.  The following is the analogue of Theorem 6.7 in \cite{BobokZN:mfold}. 
 
 \begin{lemma}\label{lem:combine}
 	Suppose \fTT{} is a \comfold{} selfmap of a tree.  
	
	Under Assumptions \ref{assume:assumption1} and \ref{assume:noncomfold}, 
	there exists a cocountable open set $\Y\subset\T$ and an \msec{} on \Y{} such that the
	associated $m$-shift system \Shiftsystem{} satisfies:
	\begin{enumerate}
		\item \fof{\kerH} is (at most) countable;
		\item every minimal set in \centOH{} is a periodic orbit;
		\item there is a set $W$ which contains every periodic orbit in \centOH{} and divides 
		\Shiftsystem{}.
	\end{enumerate}

 \end{lemma}
 
 \begin{proof}
 	By \refer{prop}{monospan} there exists a cocountable open set 
	$Y\subset\T\setminus[\Vertices\cup\fof{\Vertices}]$ 
	and a monotone, spanning \msec{} $\psi$ on \Y{} for which the corresponding 
	$m$-shift system \Shiftsystem{} satisfies 
		\begin{equation}\label{eqn:fker}
			\fof{\kerH}\cap\Y=\emptyset,
		\end{equation}
	which immediately implies the first condition
	above.  
	
	In particular, this also implies that any minimal set of $f$ which is contained in the core 
	\centOH{} (and hence
	in the kernel) is a periodic orbit, since otherwise it would have to be uncountable.  
	
	But a combination of \refer{cor}{badword} 
	and \refer{prop}{locsepset} (together with \refer{rmk}{locsepset} and \refer{lem}{unionlocdiv}) 
	gives the existence of a set
	$W$ which contains all periodic orbits in the core and which locally divides \Shiftsystem{}.
\end{proof}	
	
By \refer{rmk}{assumptions}, we can assume \refer{assume}{assumption1} holds, while failure of 
\refer{assume}{noncomfold} clearly gives the conclusion of \refer{thm}{main}.
But this means that \refer{lem}{combine} gives us the hypotheses of \refer{thm}{mshiftent}, and \refer{thm}{main} follows immediately.

\bibliographystyle{amsalpha}
\bibliography{mFold}

 \end{document}